\documentclass[a4paper]{article}
\usepackage[margin=25mm]{geometry}
\usepackage{graphicx}
\usepackage{subfigure}
\usepackage{float}
\usepackage{graphicx}
\usepackage{amssymb}
\usepackage{mathrsfs}
\usepackage{fullpage}
\usepackage{amsmath}
\usepackage{graphicx,caption}
\usepackage{eucal}
\usepackage{lipsum}
\usepackage{amsfonts}
\DeclareMathAlphabet{\mathpzc}{OT1}{pzc}{m}{it}
\usepackage{epsfig}
\usepackage{adjustbox}
\allowdisplaybreaks
\usepackage{latexsym}
\usepackage[nodisplayskipstretch]{setspace}
\makeatletter
\def\ps@pprintTitle{%
\let\@oddhead\@empty
\let\@evenhead\@empty
\def\oddfoot\@empty
\let\evenfoot\@oddfoot}
\makeatother
\usepackage{verbatim}
\providecommand{\keywords}[1]
{
  \small	
  \textbf{\textit{Keywords---}} #1
}

\title{Optimization of traffic control in $\textrm{\it MMAP[2]/PH[2]/S}$ priority queueing model with $\textrm{\it PH}$ retrial times and the preemptive repeat  policy}
\author{Raina Raj$^{1}$, Vidyottama Jain$^{1}$   \\
        \small Central University of Rajasthan, Ajmer, India$^{1}$\\
    }
\date{} % Comment this line to show today's date

\begin{document}
\maketitle

\begin{abstract}
\textbf{The presented study elaborates a multi-server  priority  queueing model considering the preemptive repeat policy and phase-type distribution ($\textit{P\!H}$) for retrial process. The incoming heterogeneous calls are categorized as handoff calls and new calls. The arrival and service processes of both types of calls follow marked Markoian arrival process ($\textit{M\!M\!A\!P}$) and  $\textit{P\!H}$ distribution with distinct parameters, respectively. An arriving new call will be blocked when all the channels are occupied, and consequently will join the orbit (virtual space) of infinite capacity to retry following $\textit{P\!H}$ distribution.  When all the channels are occupied and a handoff call arrives at the system, out of the following  two scenarios, one might take place. In the first scenario,  if all the channels are occupied with handoff calls, the arriving handoff call will be lost from the system. While in the second one, if all the channels are occupied and at least one of them is serving a new call,  the arriving handoff call will be  provided service by using preemptive priority over that new call and the preempted  new call will join the orbit.  Behaviour of the proposed system is modelled by the level dependent quasi-birth-death $\textit{(L\!D\!Q\!B\!D)}$ process. The expressions of various performance measures have been derived for the numerical illustration. An optimization problem for optimal channel allocation and traffic control has been formulated and dealt by employing appropriate heuristic approaches.}
\end{abstract} \hspace{10pt}

%TC:ignore
%\subclass{60K25, 68M20, 90B22, 90B18}

\section{Introduction} \label{section1}
In modern wireless cellular networks, priority queueing models are effectively used in various applications  where the prioritization of traffic is essential, e.g., priority of voice traffic is ensured for voice and data transmission in multiprocessor switching. Such prioritization of traffic can be classified depending  on the various features of the system, e.g., arrival discipline, service discipline, types of services, etc. 
In the literature, two types of  priority queueing models (\textit{P\!Q\!M\!s}) based on the service discipline, named as the non-preemptive priority queueing model (\textit{N\!P\!P\!Q\!M}) and the preemptive priority queueing model (\textit{P\!P\!Q\!M}),  have been discussed (refer, \cite{brandwajn2017multi,chang1965preemptive,krishnamoorthy2008map,machihara1995bridge}, etc.). In the \textit{N\!P\!P\!Q\!M}, the service of a lower priority traffic is not interrupted by the arrival of a higher priority traffic. Whereas, in the \textit{P\!P\!Q\!M}, the service of an ongoing lower priority traffic is terminated 
 by the arrival of a higher priority traffic. This higher priority traffic  gets the service in place of the preempted lower priority traffic. %(\cite{chang1965preemptive}). 
In this work, the incoming traffic is categorized  into handoff calls and new calls where  preemptive priority is assigned to 
handoff calls over new calls. The preempted new call  joins the orbit and  retries for the service.  The preemptive policy may be  distinguished as  the preemptive repeat priority  policy or as the preemptive resume priority policy. In the former one,  the service of a preempted new call is  commenced from the scratch and in the later one, the service of a preempted new call is resumed from the point at which it has been interrupted. 
%ongoing service of a new call is preempted, 
%In that scenario, the preemptive discipline must distinguish two cases: The service of a preempted new call will either commenced from the scratch or resume from the point at which it has been interrupted. The former discipline is referred as preemptive repeat priority  and the later one as preemptive resume priority.  

In this study, the preemptive repeat priority policy has been adopted for the prioritization of handoff calls over new calls.
An overview of research works on  the preemptive repeat  \textit{P\!Q\!M} can be found in the articles (refer, \cite{artalejo2001stationary,chang1965preemptive,drekic2001reducing,fiems2018exhaustive}, etc.). In these studies,  challenges arising due to the consideration of the preemptive repeat  policy have been emphasized.     However, the applicability of these models has been diminished in the present scenario, since the incoming call arrival follows Poisson process and service times is considered exponentially distributed. %(\cite{artalejo2001stationary,drekic2001reducing}). 
%breuer2002map, asmussen2017preemptive

In the cutting edge wireless technologies, the input flow of calls  possesses burstiness and correlation properties rather than the memory-less property of stationary Poisson flow. %Therefore, $\textit{M\!M\!A\!P}$ is considered as a very acknowledged stochastic process to capture the correlation among the heterogeneous calls.systems over priority queues with multi-servers,
Despite the enormous number of single/multi server \textit{P\!P\!Q\!M\!s},  the literature devoted to the analysis of such models with the consideration of more general arrival and service processes, e.g., Markovian arrival process ($\textit{M\!A\!P}$), marked Markovian arrival process ($\textit{M\!M\!A\!P}$) and phase-type ($\textit{P\!H}$) distribution, appears rather moderate. Some of the relevant studies are as follows. 	He and Alfa \cite{he1998mmap} derived stationary distribution for a single-server queueing system following $\textit{M\!M\!A\!P}$ for the arrival process and different $\textit{P\!H}$ distributions for the service of distinct classes of customers.  He and Li \cite{he2003stability} obtained stability conditions for a single-server preemptive repeat \textit{P\!Q\!M}  by assuming $\textit{M\!M\!A\!P}$ arrival and general distributed service times for distinct classes of customers. %Choi and Hwang \cite{choi1997map} also analyzed a similar single-server priority queueing model with preemptive resume priority policy for two customer classes.
%Horvath et al. \cite{horvath2009traffic} analyzed a single-server priority queueing model with $\textit{MMAP}$ arrival and  $\textit{MAP}$ for service process.
 Sun et al. \cite{sun2014map+}   presented a $\textrm{\it M\!A\!P\! +\! M\!A\!P/$M_2$/S}$  queueing model with infinite buffer considering the preemptive repeat priority policy and reservation policy. 
 %They provided absolute priority to a higher priority customer
 %by using   preemptive repeat priority policy and reservation policy. The arrival and service processes were defined with distinct parameters following $\textit{ M\!A\!P}$ and exponential distributions, respectively.
 The consideration of exponential distribution for service times restricts the applicability of the model in real life scenario.
 Recently, Klimenok et al. \cite{klimenok2020priority} proposed a  $\textrm{\it M\!M\!A\!P[2]/P\!H[2]/S}$ priority queueing model and  analyzed the prsented system  without considering the retrial behavior of the customers.

%The retrial phenomenon takes place into the system  in two cases for the new calls. The first case takes place when an arriving new call finds all the available channels occupied, and consequently joins a virtual space, called orbit, to retry  for the service after some random amount of time. These blocked new calls are referred as retrial calls in this study. And,  the second scenario occur when the service of an ongoing new call is preempted due to the arrival of handoff call, hence, it will join the orbit and will be dealt as a retrial call.but confined the model

Though, there exist a vast literature over the \textit{P\!P\!Q\!M} with single/multi server, yet a handful studies have shown the impact of retrial phenomenon over the \textit{P\!P\!Q\!M} with generalized arrival and service processes. For a detailed survey of studies over retrial queueing models, readers may refer \cite{kim2016survey} and references cited herein. Some of the relevant literature for the \textit{P\!P\!Q\!M} with retrial phenomenon is discussed here.
 Dudin et al. \cite{dudin2015priority} proposed a multi-server retrial queueing system following $\textit{M\!M\!A\!P}$ for the arrival of customers with the preemptive repeat priority policy under a random environment. They confined the  model considering exponentially distributed service and retrial times.
  %They considered exponential distributions for  service and retrial times. %They evaluated the performance of the system under a random environment.
 %, i.e., the parameters of systems will change their values if the random environment changes its state.  by applying preemptive priority policy. The arrival and service processes were denoted by $\textit{M\!M\!A\!P}$ and $\textit{P\!H}$ distribution, respectively.
   Dudin et al. \cite{dudin2016analysis}  constructed a  retrial \textit{P\!P\!Q\!M} considering $\textit{M\!M\!A\!P}$ and $\textit{P\!H}$ distribution for the arrival and service processes, respectively. Their objective was to analyze the model by considering the channel sharing scheme in cognitive radio network.
Though, the above mentioned studies considered generalized arrival and service processes yet retrial process is explained by exponential distribution only. In wireless cellular networks, the inter-retrial times are notably brief in comparison to the service times. Since, the retrial attempt is just a matter of pushing one button, these retrial customers will make numerous attempts during any given service interval. Therefore, the consideration of exponential retrial times in place of non-exponential ones could lead to under or over estimating the system parameters as shown by various studies in the literature (refer, \cite{chakravarthy2020retrial,dharmaraja2008phase,shin2011approximation}, etc.). Therefore, to obtain realistic performance measures for retrial phenomenon,  a more generalized approach, i.e., $\textrm{\it P\!H}$ distribution,  has been applied in the proposed model. %Some studies have considered $\textrm{\it P\!H}$ distributed retrial times in their models (refer, \cite{dharmaraja2008phase,shin2011approximation,chakravarthy2020retrial}, etc.), yet these studies did not consider any policy for the prioritization.

%The literature over generalized $\textrm{\it P\!H}$ distributed retrial time is restricted till studies \cite{shin2011approximation,chakravarthy2020retrial}.  

%The preempted type 2 customer has an option to join the orbit and retry for service after some random time or exit the system without obtaining the service.  

In this work,  a  $\textrm{\it M\!M\!A\!P[2]/P\!H[2]/S}$ \textit{P\!Q\!M} with   $\textrm{\it P\!H}$ distributed retrial times is introduced.  To the best of authors' knowledge, the proposed model is the first one that deals with such complex system considering the preemptive repeat priority policy. The arrival and service processes for both types of calls are described by applying $\textrm{\it M\!M\!A\!P}$ and $\textrm{\it P\!H}$ distributions with different parameters, respectively. The new call, which finds all the channels busy upon its arrival will join the orbit (virtual space) of infinite capacity and will be referred as a retrial call (\cite{jain2020numerical}). The retrial call following $\textit{P\!H}$ distribution can either retry for service or exit the system without obtaining the service. If all the channels are occupied at the arrival epoch of a handoff call, out of the following two cases, one  might occur. In the first case, the arriving handoff call will be lost from the system when all the channels are occupied with the handoff calls. In the second case, the handoff call will be provided preemptive priority over the ongoing new call when at least one of the channel is occupied with that new call. 
%In this model, an arriving handoff call is assigned priority over new call by applying preemptive repeat priority policy. As per this policy, if at the arrival epoch of a handoff call,  no idle channel is available in the system then the service of one of the new calls in service is terminated. The arriving handoff call will be considered lost from the system only if, all the channels are occupied with handoff calls.
The handoff call commenced its service in place of the preempted new call and this preempted new call  joins the orbit.  The underlying process of the system is modelled by level dependent quasi-birth-death $\textit{(L\!D\!Q\!B\!D)}$ process (\cite{bright1995calculating}). Further, a matrix analytic algorithm, proposed by Baumann and Sandmann \cite{baumann2013computing}, is applied for the analysis of the proposed model.  The detailed study over   $\textit{L\!D\!Q\!B\!D}$ process can be found in \cite{he2014fundamentals} and \cite{latouche1999introduction}. Due to the  consideration of the preemptive repeat priority policy, the dropping probability for handoff calls in the system decreases and simultaneously  the frequent termination of services for new calls increases the probability of preemption.
%decreases the dropping probability for handoff calls in the system but at the same time increases the probability of preemption for new calls also. 
Both types of probabilities (say, loss probabilities) are majorly affected by the arrival of handoff calls. Thus, an optimal value of handoff call arrival rate is estimated in order to obtain a minimum value of the total number of channels in such a way that the dropping probability and preemption probability should not exceed some pre-defined values. On the basis of the above mentioned assumptions, an optimization problem has been proposed and solved by applying appropriate heuristic algorithms.
%The incoming calls, categorized as handoff calls $(\mathcal{H})$ and new calls $(\mathcal{N})$, follow  $\textrm{\it MMAP}$. The service time distribution of handoff call and new call are described by \textit{$PH_1$} and \textit{$PH_2$}  distributions with different intensities, respectively. These new call are referred as retrial calls.The retrial call can either retry for the service after some time or exit the cell without obtaining the service. The retrial process for retrial call follows \textit{PH} distribution.The handoff calls which find atleast  one channel occupied by a new call, will be served by using preemptive repeat priority policy. As per this policy, the service of an ongoing new call will be terminated if, at the arrival epoch of a handoff call, no free channel is available. The handoff call will start service in the place of the pushed out new call. The new call whose service has been preempted will  join the orbit and retry for the service like a retrial call. If the arriving handoff call find all channels occupied with the handoff calls, then it will be lost from the system.

The layout of this work is  arranged in seven sections.  In Section \ref{section2},  a  $\textrm{\it M\!M\!A\!P[2]/P\!H[2]/S}$ model with  $\textrm{\it P\!H}$ distributed retrial times  is described.  In Section \ref{section3},  the infinitesimal generator matrix for the proposed $\textit{L\!D\!Q\!B\!D}$ process has been derived and steady-state probabilities has been computed through matrix analytic algorithm. In Section \ref{section4}, formulas of key performance measures to analyse the system efficiency  are derived explicitly. Numerical  illustrations to point out the impact of various intensities over the system performance are presented in Section \ref{section5}. An optimization problem has been formulated to evaluate the  behaviour of the system in Section \ref{section6}. Finally, the underlying model is concluded with the insight for the future works in Section \ref{section7}.

\section{Model Description} \label{section2}
This work considers a multi-server  \textit{P\!Q\!M} with the preemptive repeat policy and $\textrm{\it P\!H}$ distributed retrial process. All the other assumptions are described as follows.

\begin{figure}
	\centering
	\includegraphics[width=6.2in, height=3.1in]{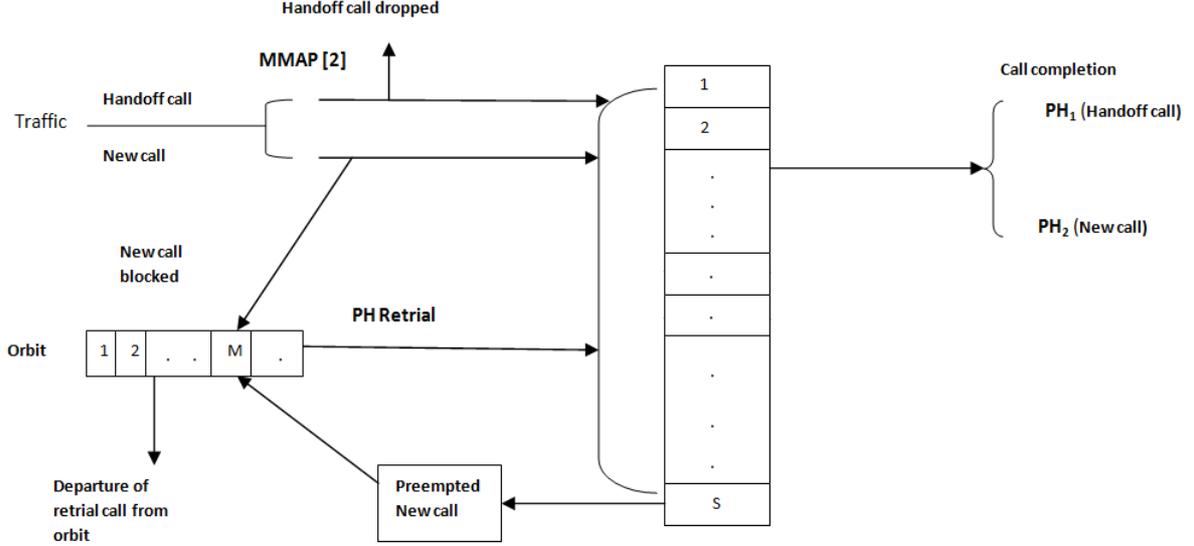} 
	\caption{A multi-server $\textrm{\it MMAP[2]/PH[2]/S}$ model  with $\textrm{\it P\!H}$ retrial times and the preemptive repeat priority policy.}
	\label{fig:fig_sim1}
\end{figure}

\begin{itemize}
	\item \textbf{Arrival Process:}
	
	\hspace{0.5cm}  The arrival of a handoff call/new call follows  \textit{M\!M\!A\!P}. The arrival in \textit{M\!M\!A\!P} is directed by the underlying process \{$\nu_t, t \geq 0$\}, which is an irreducible continuous time Markov chain with the state space of dimension $L$. The \textit{M\!M\!A\!P} is  defined by the matrices $C_0, C_{\mathcal{N}}, C_{\mathcal{H}}$; where $C_0$ represents the rates of transitions due to the occurrence of no arrival. The rates of transitions accompanied by arrival of a new call and a handoff call are described by $C_{\mathcal{N}}$ and $C_{\mathcal{H}}$, respectively.  The steady-state vector $\pi$ is the unique solution to the system $\pi C = 0,~ \pi e = 1,$ where $C = C_0 + C_{\mathcal{N}} + C_{\mathcal{H}}.$ Here  0 is a row vector consisting of 0’s of appropriate dimension and $e$ is a column vector consisting of 1’s. The fundamental arrival rates of handoff and new call are given by $\lambda_{\mathcal{H}} = \pi C_{\mathcal{H}} e$ and $\lambda_{\mathcal{N}} = \pi C_{\mathcal{N}} e$, respectively. The total fundamental arrival rate is $\lambda = \lambda_{\mathcal{H}} + \lambda_{\mathcal{N}}.$ For more details over \textit{M\!M\!A\!P}, authors suggest to refer \cite{he1996queues}.\\

	\item \textbf{Service Process:}
	
	\hspace{0.5cm} The system provides different types of services to handoff calls and new calls which follow \textrm{\it $P\!H$} distributions with distinct parameters.
	%	two types of service  for handoff call and new call with \textit{$P\!H_1$} and  \textrm{\it $P\!H_2$} distributions, respectively. 
	The service times of a handoff call following \textit{$P\!H$} distribution has $(\beta_{\mathcal{H}}, A_{\mathcal{H}})$ representation with dimension $M_{\mathcal{H}}$, i.e., $A_{\mathcal{H}}e + A_{\mathcal{H}}^0 =0.$ The average service rate of a handoff call is given by $1/\mu_{\mathcal{H}} = -\beta_{\mathcal{H}} (A_{\mathcal{H}})^{-1}e.$
	Similarly, the service times of a new call  follows \textit{$P\!H$} distribution with representation $(\beta_{\mathcal{N}}, A_{\mathcal{N}})$ of dimension $M_{\mathcal{N}}$, i.e., $A_{\mathcal{N}}e + A_{\mathcal{N}}^0 =0.$ The average service rate of a new call is given by $1/\mu_{\mathcal{N}} = -\beta_{\mathcal{N}} (A_{\mathcal{N}})^{-1}e.$\\ %and $\mu = \mu_{\mathcal{H}}+\mu_{\mathcal{N}}.$
	
	\item \textbf{Retrial Process:}

	\hspace{0.5cm}  The retrial times of a retrial call is \textit{$P\!H$} distributed with representation $(\gamma, \Gamma)$ and dimension $N$, i.e., $\Gamma e + \Gamma^0(1) + \Gamma^0(2)=0.$ Here $\Gamma^0(1)$ shows the absorption due to the departure from the cell and $\Gamma^0(2)$ denotes the absorption due to the retrial attempt. The average retrial rate is given by $1/ \theta = -\gamma (\Gamma)^{-1} e. $
	
\end{itemize}

\section{Mathematical Analysis} \label{section3}

% i = \mathpzc{l}; j = \kappa; k= \mathfrak{j}; l = \mathfrak{i} 

The underlying process \{$\Xi(t), t \geq 0 \}$ for a cell is defined by the following state space:
\[ \Omega = \{(\mathpzc{l}, \kappa,   \mathfrak{j}, v, s_{\mathcal{H}}, s_{\mathcal{N}}, \mathpzc{r}); \mathpzc{l} \geq 0,~0 \leq \kappa \leq S,~ 0 \leq \mathfrak{j} \leq S,~1 \leq v \leq L\},\]
where, 
\begin{itemize}
	\item $\mathpzc{l}$ is the number of retrial calls,
	\item $\kappa$ is the  number of busy channels,
	\item $\mathfrak{j}$ is the number of handoff calls in the system receiving service, 
	\item $v$ is the current phase of $\textrm{\it M\!M\!A\!P}$,
	\item $s_{\mathcal{H}}= \{(s_{\mathcal{H}}^1, s_{\mathcal{H}}^2,\ldots,s_{\mathcal{H}}^{\mathfrak{j}}); 1 \leq s_{\mathcal{H}}^{\nu_1} \leq M_\mathcal{H}; 1 \leq \nu_1 \leq \mathfrak{j}\},$ where the $\nu_1^{\textrm{th}}$  handoff call,  out of  $\mathfrak{j}$ number of  handoff calls, is being served in phase $s_{\mathcal{H}}^{\nu_1},$
	\item $s_{\mathcal{N}}= \{(s_{\mathcal{N}}^1, s_{\mathcal{N}}^2,\ldots,s_{\mathcal{N}}^{\kappa-\mathfrak{j}}); 1 \leq s_{\mathcal{H}}^{\nu_2} \leq M_\mathcal{N}; 1 \leq \nu_2 \leq \kappa-\mathfrak{j}\},$ where the $\nu_2^{\textrm{th}}$  new call,  out of  $\kappa-\mathfrak{j}$ number of new calls, is being served in phase $s_{\mathcal{N}}^{\nu_2},$
	\item $\mathpzc{r} = \{(r^1, r^2,\ldots, r^{\mathpzc{l}}); 1 \leq r^{h} \leq N; 1 \leq h \leq \mathpzc{l}\},$ where the $h^{\textrm{th}}$ retrial call  in  the orbit is in phase $r^h.$

\end{itemize}
The stochastic process \{$\Xi(t), t \geq 0 \}$ can be modelled as  $\textrm{\it L\!D\!Q\!B\!D}$ process with the tri-diagonal infinitesimal generator matrix provided as follows:
\begin{center}
	$\mathscr{Q} =
	\begin{pmatrix}
	\mathscr{Q}_{0,0} & \mathscr{Q}_{0,1} & 0 & 0 & 0 &  \\
	\mathscr{Q}_{1,0}& \mathscr{Q}_{1,1}& \mathscr{Q}_{1,2} & 0 & 0 &  \\
	0  & \mathscr{Q}_{2,1} & \mathscr{Q}_{2,2} & \mathscr{Q}_{2,3} & 0 &   \\
%	0 & 0   & \mathscr{Q}_{3,2} & \mathscr{Q}_{3,3} &\mathscr{Q}_{3,4} &               \\
	%0 & 0 & 0  &  &  &               \\
	%&  &     &          & \\
	% &   &   \ddots & \ddots & \ddots & \\
	%      &  &   Q_{i,i-1}  &  Q_{i,i} & Q_{i,i+1} &  \\ 
	%     &  & &  Q_{i,i-1}  &  Q_{i,i} & Q_{i,i+1}  \\   
	&   &  & \ddots & \ddots & \ddots \\
	&   & & & \ddots & \ddots & \ddots 
	\end{pmatrix}.$\\
\end{center}

We define the following notations in order to carry out the analysis. 
\begin{itemize}
	\item $I_u$ is an identity matrix of dimension $u$.
	\item ${\bf O}$ is a  zero matrix of appropriate dimension.
	\item O(A) denotes the order of matrix A.
%	\item $I(u,z) = \underbrace{I_u \otimes I_u \otimes....\otimes I_u}_z.$
%	\item $0(u,z) = \underbrace{0_u \oplus 0_u \oplus \ldots \oplus 0_u}_z.$
	\item $\Psi_{\mathcal{N}}(\kappa-\mathfrak{j}) = \underbrace{A_{\mathcal{N}} \oplus A_{\mathcal{N}} \oplus \ldots \oplus A_{\mathcal{N}}}_{\kappa-\mathfrak{j}}$ represents that $\kappa-\mathfrak{j}$ number of new calls are receiving services.
	
	\item  $\Psi_{\mathcal{H}}(\mathfrak{j}) = \underbrace{A_{\mathcal{H}} \oplus A_{\mathcal{H}} \oplus \ldots \oplus A_{\mathcal{H}}}_{\mathfrak{j}}$ represents that $\mathfrak{j}$ number of handoff calls are receiving services.
	
	\item $ \Psi_{orbit}(\mathpzc{l}) = \underbrace{\Gamma \oplus \Gamma \oplus \ldots \oplus \Gamma}_{\mathpzc{l}}$ 
	represents that $\mathpzc{l}$ number of  new calls are  retrying.
	\item $\hat \Psi_{orbit}(\mathpzc{l}) = \underbrace{(\Gamma^0(2)\gamma) \oplus (\Gamma^0(2)\gamma) \oplus \ldots \oplus (\Gamma^0(2)\gamma)}_{\mathpzc{l}}$ represents that $\mathpzc{l}$ number of new calls are having unsuccessful retrial attempts.  
	\item $\Phi_{\mathcal{N}}(\kappa-\mathfrak{j}) = \displaystyle{\sum_{y=0}^{\kappa-\mathfrak{j}-1}I_{M_{\mathcal{N}}^y} \otimes A_{\mathcal{N}}^0 \otimes I_{M_{\mathcal{N}}^{ \kappa-\mathfrak{j}-y-1}}}$ represents  that any one out of the $\kappa-\mathfrak{j}$ number of new calls has completed the service.
	
	\item $\displaystyle{\Phi_{\mathcal{H}}(\mathfrak{j}) = \sum_{y=0}^{\mathfrak{j}-1}I_{M_{\mathcal{H}}^y} \otimes A_{\mathcal{H}}^0 \otimes I_{M_{\mathcal{H}}^{ \mathfrak{j}-y-1}}}$ represents that  any one out of the $\mathfrak{j}$ number of handoff calls has completed the service.
	
	\item $\Phi_{orbit}(\mathpzc{l}+1) = \displaystyle{\sum_{y=0}^{\mathpzc{l}}I_{N^y} \otimes \Gamma^0(1) \otimes I_{N^ {\mathpzc{l}-y}}}$ represents that any one out of the $\mathpzc{l}$ number of new calls leaves the orbit as well as the cell without getting connected.
	
	\item $\hat \Phi_{orbit}(\mathpzc{l}+1) = \displaystyle{\sum_{y=0}^{\mathpzc{l}}I_{N^y} \otimes (\Gamma^0(2) \otimes \beta_{\mathcal{N}}) \otimes I_{N^{ \mathpzc{l}-y}}}$ represents that any one out of the $\mathpzc{l}$ number of new calls is getting service after its  successful retrial attempt. 
\end{itemize}

	Note that notations $\oplus$ and  $\otimes$ are used for the Kronecker sum and  the   Kronecker product of two matrices, respectively. For more description over Kronecker sum and Kronecker product, refer  \cite{dayar2012analyzing}.
	
The intensities of the upper diagonal of the $\mathscr{Q}$ matrix represent the scenario when
 one new  call joins the orbit due to the non availability of idle channels. The intensities of the lower diagonal show the loss of one retrial call either due to the successful retrial or due to the departure from the orbit without obtaining the service. The main diagonal represents transitions due to the arrival or service of handoff calls and new calls or the transitions of a retrial call from one phase to another phase. The number of retrial calls is  not  changed during these transitions.  The block matrices are defined as follows.
\begin{align*}
& \textbf {Upper  Diagonal :}\\
&\mathscr{Q}_{\mathpzc{l},\mathpzc{l}+1} = \text{diag}\{X_\mathpzc{l}(0), X_\mathpzc{l}(1), \ldots, X_\mathpzc{l}(S)  \};~~\mathpzc{l} \geq 0,\\
%	& \text{where}\\
	& X_\mathpzc{l}(\kappa) = {\bf O}; ~~ \kappa = \overline{0,S-1},\\
%	& \text{${\bf O}$ is a  zero matrix of appropriate dimension,}\\
	& X_\mathpzc{l}(S) = 
	\begin{pmatrix}
		X_\mathpzc{l}(S,0) & \hat{X}_\mathpzc{l}(S,0) & 0   & \cdots &0&0  \\
		0& X_\mathpzc{l}(S,1)& \hat{X}_\mathpzc{l}(S,1) & \cdots  & 0&0 \\
	%	0  & 0 & X_\mathpzc{l}(\kappa,2)  & \cdots &  0&0 \\
		\vdots      & \vdots & \vdots & \ddots  & \vdots& \vdots\\
		0  & 0 & 0 & \cdots& X_\mathpzc{l}(S,S-1) &\hat{X}_\mathpzc{l}(S,S-1)\\
		0    &  0    &  0&     \hdots  & 0 & X_\mathpzc{l}(S,S)
	\end{pmatrix},\\
	& X_\mathpzc{l}(S,\mathfrak{j}) = (C_{\mathcal{N}} \otimes I_{M_{\mathcal{H}}^{\mathfrak{j}}M_{\mathcal{N}}^{S-\mathfrak{j}} N^{\mathpzc{l}}}) \otimes \gamma ;~ \forall \mathpzc{l} \geq 0,~~\mathfrak{j} = \overline{0,S},\\
	& \hat{X}_\mathpzc{l}(S,\mathfrak{j}) = (C_{\mathcal{H}} \otimes I_{M_{\mathcal{H}}^{\mathfrak{j}}}  )\otimes \beta_{H}\otimes I_{M_{\mathcal{N}}^{S-\mathfrak{j}-1}} \otimes e(M_{\mathcal{N}})\otimes  I_{N^{\mathpzc{l}}}\otimes \gamma ;~ \forall \mathpzc{l} \geq 0,~~\mathfrak{j} = \overline{0,S-1}.\\
&	 \textbf  {Lower Diagonal :}\\
 &\mathscr{Q}_{\mathpzc{l}+1,\mathpzc{l}} =
	\begin{pmatrix}
		Z_\mathpzc{l}(0) & \hat{Z}_\mathpzc{l}(0) & 0   & \cdots &0&0  \\
		0& Z_\mathpzc{l}(1)& \hat{Z}_\mathpzc{l}(1) & \cdots  & 0&0 \\
%		0  & 0 & Z_\mathpzc{l}(2)  & \cdots &  0&0 \\
		\vdots      & \vdots & \vdots & \ddots  & \vdots& \vdots\\
		0  & 0 & 0 & \cdots& Z_\mathpzc{l}(S-1) &\hat{Z}_\mathpzc{l}(S-1)\\
		0    &  0    &  0&     \hdots  & 0 & Z_\mathpzc{l}(S)
	\end{pmatrix}; ~~\mathpzc{l} \geq 0,\\
%	& \text{where}\\
&Z_\mathpzc{l}(\kappa) = \text{diag}\{	Z_\mathpzc{l}(\kappa,0), Z_\mathpzc{l}(\kappa,1), \ldots, Z_\mathpzc{l}(\kappa,\kappa) \};~~\kappa = \overline{0,S},~ O(Z_\mathpzc{l}(\kappa)) = (\kappa+1)\times (\kappa+1),\\	 %\text{$Z_\mathpzc{l}(\kappa)$ is a square matrix of order $(\kappa+1)\times (\kappa+1)$,}\\
& Z_\mathpzc{l}(\kappa,\mathfrak{j}) = I_{LM_{\mathcal{H}}^{\mathfrak{j}}M_{\mathcal{N}}^{\kappa-\mathfrak{j}}} \otimes \Phi_{orbit}(\mathpzc{l}+1);~ \forall \mathpzc{l} \geq 0,~  \kappa = \overline{0,S},~ \mathfrak{j} = \overline{0,\kappa},\\
& \hat{Z}_\mathpzc{l}(\kappa) = \text{diag}\{	\hat{Z}_\mathpzc{l}(\kappa,0), \hat{Z}_\mathpzc{l}(\kappa,1), \ldots, \hat{Z}_\mathpzc{l}(\kappa,\kappa-1) \};~~\kappa = \overline{0,S-1},~ O(\hat{Z}_\mathpzc{l}(\kappa)) = (\kappa+1)\times (\kappa+2),\\
	& \hat{Z}_\mathpzc{l}(\kappa,\mathfrak{j}) = I_{LM_{\mathcal{H}}^{\mathfrak{j}}M_{\mathcal{N}}^{\kappa-\mathfrak{j}}} \otimes \hat{\Phi}_{orbit}(\mathpzc{l}+1);~ \forall \mathpzc{l} \geq 0,~  \kappa = \overline{0,S-1},~\mathfrak{j} = \overline{0,\kappa}.\\
	   & \textbf  {Main Diagonal :}\\
	&\mathscr{Q}_{\mathpzc{l},\mathpzc{l}} =
	\begin{pmatrix}
		Y_\mathpzc{l}(0) & \hat{Y}_\mathpzc{l}(0) & 0  & \cdots &0&0  \\
		\overline{Y}_\mathpzc{l}(1)& Y_\mathpzc{l}(1)& \hat{Y}_\mathpzc{l}(1)& \cdots  & 0&0 \\
	%	0  & \overline{Y}_\mathpzc{l}(2) & Y_\mathpzc{l}(2)  &  \hat{Y}_\mathpzc{l}(2)&\cdots &  0&0 \\
		\vdots      & \vdots & \vdots & \ddots  & \vdots& \vdots\\
		0  & 0 & 0 & \cdots& Y_\mathpzc{l}(S-1) &\hat{Y}_\mathpzc{l}(S-1)\\
		0    &  0 & 0   &       \hdots  & \overline{Y}_\mathpzc{l}(S) & Y_\mathpzc{l}(S)
	\end{pmatrix}; ~~\mathpzc{l} \geq 0,\\
%& where\\
	& \hat{Y}_\mathpzc{l}(\kappa) =    \begin{pmatrix}
		\hat{Y}_\mathpzc{l}(\kappa,0,\mathcal{N}) & \hat{Y}_\mathpzc{l}(\kappa,0,\mathcal{H}) & 0   & \cdots &0&0  \\
		0& \hat{Y}_\mathpzc{l}(\kappa,1,\mathcal{N})& \hat{Y}_\mathpzc{l}(\kappa,1,\mathcal{H}) & \cdots  & 0&0 \\
	%	0  & 0 & \hat{Y}_\mathpzc{l}(\kappa,2,\mathcal{N})  & \cdots &  0&0 \\
		\vdots      & \vdots & \vdots & \ddots  & \vdots& \vdots\\
		0  & 0 & 0 & \cdots& \hat{Y}_\mathpzc{l}(\kappa,\kappa,\mathcal{N}) &\hat{Y}_\mathpzc{l}(\kappa,\kappa,\mathcal{H})
	\end{pmatrix};~~\kappa = \overline{0,S-1},\\
	&  O(\hat{Y}_\mathpzc{l}(\kappa))=(\kappa+1)\times (\kappa+2),\\
%	& \text{$\hat{Y}_\mathpzc{l}(\kappa)$ is a rectangular matrix of order $(\kappa+1)\times (\kappa+2)$,}\\
	& \hat{Y}_\mathpzc{l}(\kappa,\mathfrak{j},\mathcal{N}) = ((C_{\mathcal{N}} \otimes I_{M_{\mathcal{H}}^{\mathfrak{j}}M_{\mathcal{N}}^{\kappa-\mathfrak{j}}})\otimes \beta_{\mathcal{N}})\otimes I_{N^{\mathpzc{l}}};~\forall \mathpzc{l} \geq 0,~  \kappa = \overline{0,S-1},~ \mathfrak{j} = \overline{0,\kappa},\\
	&\hat{Y}_\mathpzc{l}(\kappa,\mathfrak{j},\mathcal{H}) = (C_{\mathcal{H}} \otimes I_{M_{\mathcal{H}}^{\mathfrak{j}}}\otimes \beta_{\mathcal{H}}) \otimes I_{M_{\mathcal{N}}^{\kappa-\mathfrak{j}} N^{\mathpzc{l}}};~ \forall \mathpzc{l} \geq 0,~  \kappa = \overline{0,S-1},~ \mathfrak{j} = \overline{0,\kappa},\\
& Y_\mathpzc{l}(\kappa) = \text{diag}\{	Y_\mathpzc{l}(\kappa,0), Y_\mathpzc{l}(\kappa,1), \ldots, Y_\mathpzc{l}(\kappa,\kappa) \};~~\kappa = \overline{0,S},~ O(Y_\mathpzc{l}(\kappa)) = (\kappa+1)\times (\kappa+1),\\
	& Y_\mathpzc{l}(\kappa,\mathfrak{j}) = C_0 \oplus \Psi_{\mathcal{H}}(\mathfrak{j})\oplus \Psi_{\mathcal{N}}(\kappa-\mathfrak{j})\oplus \Psi_{orbit}({\mathpzc{l}});~ \forall \mathpzc{l} \geq 0,~  \kappa = \overline{0,S-1},~  \mathfrak{j} = \overline{0,\kappa},\\
	& Y_\mathpzc{l}(S,\mathfrak{j}) = C_0 \oplus\Psi_{\mathcal{H}}(\mathfrak{j})\oplus \Psi_{\mathcal{N}}(S-\mathfrak{j})\oplus \Psi_{orbit}({\mathpzc{l}}) + I_{LM_{\mathcal{H}}^{\mathfrak{j}}M_{\mathcal{N}}^{S-\mathfrak{j}}} \otimes \hat{\Psi}_{orbit}({\mathpzc{l}}) ;~ \forall \mathpzc{l} \geq 0,~  \mathfrak{j} = \overline{0,S-1},\\
	& Y_\mathpzc{l}(S,S) = (C_0+C_\mathcal{H}) \oplus\Psi_{\mathcal{H}}(S)\oplus  \Psi_{orbit}({\mathpzc{l}}) + I_{LM_{\mathcal{H}}^{S}} \otimes \hat{\Psi}_{orbit}({\mathpzc{l}}) ;~ \forall \mathpzc{l} \geq 0,\\
	%	& and \\
	& \overline{Y}_\mathpzc{l}(\kappa)  = 
	\begin{pmatrix}
		\overline{Y}_\mathpzc{l}(\kappa,0,\mathcal{N})  & 0 & \cdots & 0& 0 \\
		\overline{Y}_\mathpzc{l}(\kappa,1,\mathcal{H})  & \overline{Y}_\mathpzc{l}(\kappa,1,\mathcal{N}) & \cdots & 0&0  \\
		\vdots &   \vdots  & \ddots & \vdots&\vdots\\
		0  & 0 & \cdots&   0& \overline{Y}_\mathpzc{l}(\kappa,\kappa-1,\mathcal{N})\\
		0  & 0 &   \cdots&  0&\overline{Y}_\mathpzc{l}(\kappa,\kappa,\mathcal{H})
	\end{pmatrix};~~\kappa = \overline{1,S},\\
	&O(\overline{Y}_\mathpzc{l}(\kappa))=(\kappa+1)\times \kappa,  \\
%	& \text{$\overline{Y}_\mathpzc{l}(\kappa)$ is a rectangular matrix of order $(\kappa+1)\times \kappa$,}\\
	&  \overline{Y}_\mathpzc{l}(\kappa,\mathfrak{j},\mathcal{N}) = I_{LM_{\mathcal{H}}^{\mathfrak{j}}}\otimes  \Phi_{\mathcal{N}}(\kappa-\mathfrak{j})  \otimes I_{N^{\mathpzc{l}}};~ \forall \mathpzc{l} \geq 0,~  \kappa = \overline{1,S},~   \mathfrak{j} = \overline{0,\kappa-1},\\
	&  \overline{Y}_\mathpzc{l}(\kappa,\mathfrak{j},\mathcal{H}) = I_{L} \otimes \Phi_{\mathcal{H}}(\mathfrak{j}) \otimes I_{M_{\mathcal{N}}^{\kappa-\mathfrak{j}} N^{\mathpzc{l}}};~\forall \mathpzc{l} \geq 0,~  \kappa = \overline{1,S},~   \mathfrak{j} = \overline{1,\kappa}.
\end{align*}
Since the structure of the generator matrix is complex, a compact analytical form of the steady-state distributions is strenuous to achieve. Therefore, an algorithmic procedure is adopted to solve the proposed $\textrm{\it LDQBD}$ process.

%	It can be observed that  a closed-form analytical solution  of the steady-state  distribution is intractable for the $\textrm{\it LDQBD}$ process. We now describe the algorithmic procedure adopted here.

	\subsection{Steady-State Analysis}
Let $\displaystyle{{z_s} = \{{z_s}(0), {z_s}(1),{z_s}(2),\ldots, {z_s}(M-1), {z_s}(M),\ldots,\ldots\}}$ be the steady-state probability vector of generator matrix $\mathscr{Q}$ satisfying $\displaystyle{{z_s} \mathscr{Q} = 0; {z_s} e =1.}$ Here, element ${z_s}(\mathpzc{l})$ contains\\ $\displaystyle{1 \times  \sum_{\kappa=0}^{S}\sum_{\mathfrak{j}=0}^{\kappa} LM_{\mathcal{H}}^{\mathfrak{j}} M_{\mathcal{N}}^{\kappa-\mathfrak{j}}N^{\mathpzc{l}}}$ vector components,   $\displaystyle{{z_s}(\mathpzc{l}) = \{{z_s}(\mathpzc{l},0),~ {z_s}(\mathpzc{l},1),~ {z_s}(\mathpzc{l},2), \ldots, {z_s}(\mathpzc{l},\kappa)\}}$; ${z_s}(\mathpzc{l},\kappa)$ consists of $\displaystyle{1 \times \sum_{\mathfrak{j}=0}^{\kappa}  LM_{\mathcal{H}}^{\mathfrak{j}} M_{\mathcal{N}}^{\kappa-\mathfrak{j}}N^{\mathpzc{l}}}$ vector components,  ${z_s}(\mathpzc{l},\kappa) = \{{z_s}(\mathpzc{l}, \kappa,0),{z_s}(\mathpzc{l}, \kappa,1),{z_s}(\mathpzc{l}, \kappa,2),\ldots, {z_s}(\mathpzc{l},\kappa,\mathfrak{j})\};$  ${z_s}(\mathpzc{l},\kappa,\mathfrak{j})$ is $\displaystyle{1 \times   LM_{\mathcal{H}}^{\mathfrak{j}} M_{\mathcal{N}}^{\kappa-\mathfrak{j}}N^{\mathpzc{l}}}$ vector; $ \mathpzc{l}\geq 0,~ 0\leq \kappa \leq S,~ 0\leq \mathfrak{j} \leq \kappa.$ 

According to  Neuts \cite{neuts1994matrix},  $z_s$ satisfies the matrix-geometric relationship ${z_s}(\mathpzc{l}+1) = {z_s}(\mathpzc{l}) \Re^{(\mathpzc{l})};~~\mathpzc{l} \geq 0$, where the family of the matrices $\{\Re^{(\mathpzc{l})},~ \mathpzc{l} \geq 0\}$, called rate matrices,  are the minimal non-negative solutions to the following system of equations
\begin{align}
	\mathscr{Q}_{\mathpzc{l},\mathpzc{l}+1} + \Re^{(\mathpzc{l})} \mathscr{Q}_{\mathpzc{l},\mathpzc{l}} + \Re^{(\mathpzc{l})}(\Re^{(\mathpzc{l}+1)} \mathscr{Q}_{\mathpzc{l}+1,\mathpzc{l}}) &= 0;~ \mathpzc{l} \geq 0,\label{eq:a2}
\end{align}
and $z_s(0)$ is the solution of the equation ${z_s}(0)(\mathscr{Q}_{0,0} + \Re^{(0)} \mathscr{Q}_{1,0})=0.$ For the numerical computation, the infinite generator matrix $\mathscr{Q}$ is truncated up to a large finite number, say $M$. Once the truncation level $M$ is computed, the original system with infinite orbit size is approximated by the system with orbit size $M$. Therefore, the system will have a unique stationary distribution. Further, a matrix analytic  algorithm, provided by Baumann and Sandmann \cite{baumann2013computing}, has been applied to compute the steady-state probability vector of the proposed $\textrm{\it L\!D\!Q\!B\!D}$ system. Steps of the algorithm for computing $z_s$ will be as follows:%The algorithm for computing $z_s$ works as follows:
%Once a truncation level  $M$ is determined, the original system $\mathscr{Q}$ is approximated by the system which has  orbit of size  $M$.
%\textbf{Algorithm:} 
\begin{itemize}
\item Step 1: Define $\Re^{(M)} = 0$. Choose $M$ a large finite number such that $||\Re^{(M)}-\Re^{(M-1)} ||_{\infty} \leq \epsilon,$ where $\epsilon$ is a pre-defined positive value.
	\item Step 2: For $\mathpzc{l}= M, M-1,\ldots,1,$ compute and store
\begin{center}
    	$ \Re^{(\mathpzc{l}-1)}  = - \mathscr{Q}_{\mathpzc{l}-1,\mathpzc{l}}(\mathscr{Q}_{\mathpzc{l},\mathpzc{l}}+ \Re^{(\mathpzc{l})} \mathscr{Q}_{\mathpzc{l}+1,\mathpzc{l}})^{-1}.$\end{center}
	\item Step 3: Determine a solution of ${x_s}(0)(\mathscr{Q}_{0,0} + \Re^{(0)} \mathscr{Q}_{1,0})=0.$ 
	\item Step 4: For $\mathpzc{l}= M-1, M-2,\ldots,1,0,$ compute ${x_s}(\mathpzc{l}+1) = {x_s}(\mathpzc{l}) \Re^{(\mathpzc{l})}.$
	\item Step 5: By normalizing $x_s = (x_s(0),x_s(1),\ldots,x_s(M)),$ determine $z_s$, i.e.,
	\begin{align*}
		z_s = \frac{x_s}{c}, ~~~ where~~ c = \sum_{\mathpzc{l}=0}^{M}||x_s||,
	\end{align*}
	where $||.||$ is the row sum norm.
\end{itemize}
%The stopping criteria to get the cut-off value $M$ for the orbit capacity is $||\Re^{(\mathpzc{l})}-\Re^{(\mathpzc{l}-1)} ||_{\infty} \leq \epsilon,$ where $\epsilon$ is a pre-defined positive value.
\section{Performance Measures} \label{section4}
The following relevant  performance measures for  the proposed system  are calculated, after computing the  steady-state distribution $z_s$.
\begin{enumerate}
	\item The dropping probability of a handoff call:
	\[ P_d =  \frac{1}{\lambda_{\mathcal{H}}} \Big(\sum_{\mathpzc{l}=0}^{M-1} {z_s} (\mathpzc{l},S,S)\times C_{\mathcal{H}} e\Big). \]
	
	\item The blocking probability of a new call:
	\[ P_b =  \frac{1}{\lambda_{\mathcal{N}}} \Big(\sum_{\mathfrak{j}=0}^{S} {z_s} (M-1,S,\mathfrak{j})\times C_{\mathcal{N}} e\Big). \]
	
	\item The probability that there are $\mathfrak{j}$ number of handoff calls receiving service:
	\[P_{\mathcal{H}}(\mathfrak{j}) = \sum_{\mathpzc{l}=0}^{M-1} \sum_{\kappa=1}^{S}{z_s} (\mathpzc{l},\kappa,\mathfrak{j})e. \]
	
	\item The probability that there are $\mathfrak{j}^{'}$ number of new calls receiving service:
	\[P_{\mathcal{N}}(\mathfrak{j}^{'}) = \sum_{\mathpzc{l}=0}^{M-1} \sum_{\kappa=1}^{S}{z_s} (\mathpzc{l},\kappa,\kappa-\mathfrak{j}^{'})e. \] 
	
	\item Expected number of handoff calls receiving service:
	\[E_{\mathcal{H}} = \sum_{\mathfrak{j}=1}^{S} \mathfrak{j} P_{\mathcal{H}}(\mathfrak{j})e. \]

	\item Expected number of new calls receiving service:
	
	\[E_{\mathcal{N}} = \sum_{\mathfrak{j}^{'}=1}^{S} \mathfrak{j}^{'} P_{\mathcal{H}}(\mathfrak{j}^{'})e. \]
	
	\item The probability that there are $\mathpzc{l}$ number of retrial calls:
	\[P_{orbit}(\mathpzc{l}) = \sum_{\kappa=0}^{S}\sum_{\mathfrak{j}=0}^{\kappa} {z_s} (\mathpzc{l},\kappa,\mathfrak{j})e. \]
	
	\item Expected number of retrial calls:
	\[E_{orbit} = \sum_{\mathpzc{l}=0}^{M-1} \mathpzc{l} P_{orbit}(\mathpzc{l})e. \]
	
	\item The intensity at which both types of calls are served successfully:
	\[T_P = \sum_{\mathpzc{l}=0}^{M-1}\sum_{\kappa=1}^{S}\sum_{\mathfrak{j}=0}^{\kappa} \mu {z_s} (\mathpzc{l},\kappa,\mathfrak{j})  e. \]         
	
	\item The probability that an arriving handoff call preempts the service of an ongoing new call:
	\[P_{preempt} = \frac{1}{\lambda_{\mathcal{H}}} \sum_{\mathpzc{l}=0}^{M-1}\sum_{\mathfrak{j}=0}^{\kappa-1}  {z_s} (\mathpzc{l},S,\mathfrak{j}) \times C_{\mathcal{H}}  e. \]       
	
	\item The intensity by which a retrial call is successfully connected to an available channel:
	\[ \theta_r^{succ} =  \sum_{\mathpzc{l}=1}^{M-1}\sum_{\kappa=0}^{S-1}\sum_{\mathfrak{j}=0;\mathfrak{j} \leq \kappa}^{\kappa-1} \theta {z_s} (\mathpzc{l},\kappa,\mathfrak{j})  (e
	{\scriptstyle{(LM_{\mathcal{H}}^{\mathfrak{j}}M_{\mathcal{N}}^{\kappa-\mathfrak{j}}N^{\mathpzc{l}-1}})} \otimes( \Gamma^{0}(2)\otimes \beta_{\mathcal{N}}))e.\]
	
	\item The probability that a retrial call will exit the system without obtaining the service:
	
	\[P_{leave}^{no-service} = \frac{1}{\theta}  \Big(\sum_{\mathpzc{l}=1}^{M-1}\sum_{\kappa=0}^{S}\sum_{\mathfrak{j}=0}^{\kappa} {z_s} (\mathpzc{l},\kappa,\mathfrak{j}) (e
	{\scriptstyle{(LM_{\mathcal{H}}^{\mathfrak{j}}M_{\mathcal{N}}^{\kappa-\mathfrak{j}}N^{\mathpzc{l}-1}})} 
	\otimes \Gamma^0(1)) e\Big). \]
	
	\item The probability that an arriving new call will directly join the orbit:
	\[P_{orbit}^{join} = \frac{1}{\lambda_{\mathcal{N}}}  \Big(\sum_{\mathpzc{l}=0}^{M-1}\sum_{\mathfrak{j}=0}^{\kappa} {z_s} (\mathpzc{l},S,\mathfrak{j})  \times C_{\mathcal{N}} e\Big). \]
\end{enumerate}

In the next section,  the behaviour of the key performance measures will be explored with respect to the several rates for the presented model.

%After obtaining the explicit expressions for the performance measures, next, an illustrative behaviour of some of the performance measures will be  described in the  Section \ref{section5} for the sensitivity analysis propose. We will explore the impact of various intensities over the presented model.Due to the assumption of the $\textrm{\it PH}$ distribution of inter-retrial times, computational problems inevitably arise for larger values of $M$. Therefore,  the  proposed system is truncated up to a finite number $M$ which is computed by considering the  tolerance value $\epsilon = 10^{-5}$.

%truncation approach has been applied to compute $M$ by considering the  tolerance value $\epsilon = 10^{-5}$. %Though, the computation has been carried out for the sparse block matrices, still the calculation for the larger values of $M$ is quite difficult. %Therefore, a suitable and efficient algorithm will be developed in the future for the computation of such types of complex models.

\section{Numerical Illustration} \label{section5}
In this section, the qualitative behaviour of the proposed model is explored through a few experiments. It can be noticed that in the proposed model computational complexities inevitably arise for a large value of $M$ due to the consideration of $\textrm{\it P\!H}$ distributed inter-retrial times.
 Therefore, it is needed that the proposed model should be truncated up to a finite number $M$. For this purpose, consider $\epsilon = 10^{-5}$.
 For the numerical computation,  the matrices for the $\textrm{\it M\!M\!A\!P}$ are computed as follows,
\begin{align}
	\nonumber
	C_{0}= \begin{pmatrix}
		-1.3431 & 0.0230\\
		0 & -17.183
	\end{pmatrix},~~~C_{\mathcal{H}}= \begin{pmatrix}
		0.6600 & 0\\
		0.2567 & 8.3351
	\end{pmatrix},~~~ C_{\mathcal{N}}= \begin{pmatrix}
		0.6600 & 0\\
		0.2567 & 8.3351
	\end{pmatrix}.
\end{align}

The  correlation coefficients for both types of calls are $C_{r}^{(1)}=C_{r}^{(2)}=0.2211$ and the  variation coefficients  for both types of calls are $C_{r}^{(1)}=C_{r}^{(2)}=0.2211.$ 
The arrival rates  $\lambda_{\mathcal{H}}= \lambda_{\mathcal{N}} =1$.
Let $\textrm{\it P\!H}$ distributions parameters for the service rates of a handoff and a new call  be
\begin{align}
	\nonumber
	\beta_{\mathcal{H}}= \begin{pmatrix}
		1, &0
	\end{pmatrix},~~ A_{\mathcal{H}}= \begin{pmatrix}
		-1 & 1\\
		0 & -1
	\end{pmatrix},~~\textrm{and}~~~\beta_{\mathcal{N}}= \begin{pmatrix}
		1,& 0
	\end{pmatrix}, ~~ A_{\mathcal{N}}= \begin{pmatrix}
		-1 & 1\\
		0 & -1
	\end{pmatrix}.
\end{align}
The fundamental service rates are  $\mu_{\mathcal{H}}=\mu_{\mathcal{N}}=0.5$. The retrial rate of a retrial call, following  $\textrm{\it PH}$ distribution, is given by the parameters
\begin{align*}
	\nonumber
	\gamma= \begin{pmatrix}
		0.5, & 0.5
	\end{pmatrix},~~~ \Gamma = \begin{pmatrix}
		-2 & 2\\
		0 & -2
	\end{pmatrix},~~~~\theta=1.33.
\end{align*}
%The average retrial intensity $(\theta)$ is  1.33.\\

To demonstrate the feasibility of the developed model,
some interesting observations of the proposed system are described through the following numerical experiments. These experiments will present the behaviour of  performance measures with respect to  arrival, service and retrial rates.\\
\textbf{Experiment 1:} The objective here is to analyze the impact of arrival rate of handoff call ($\lambda_{\mathcal{H}}$) and arrival rate of new call ($\lambda_{\mathcal{N}}$) over the loss probabilities, i.e., dropping probability ($P_d$) and preemption probability ($P_{preempt}$).% associated with  the model under study.

\begin{figure}
	\centering
	\subfigure[$P_d$ versus  $\lambda_{\mathcal{H}}$]
	{\includegraphics[trim= 3cm 0.1cm 3.8cm 0.4cm, height = 5.55cm,width = 0.5\textwidth]{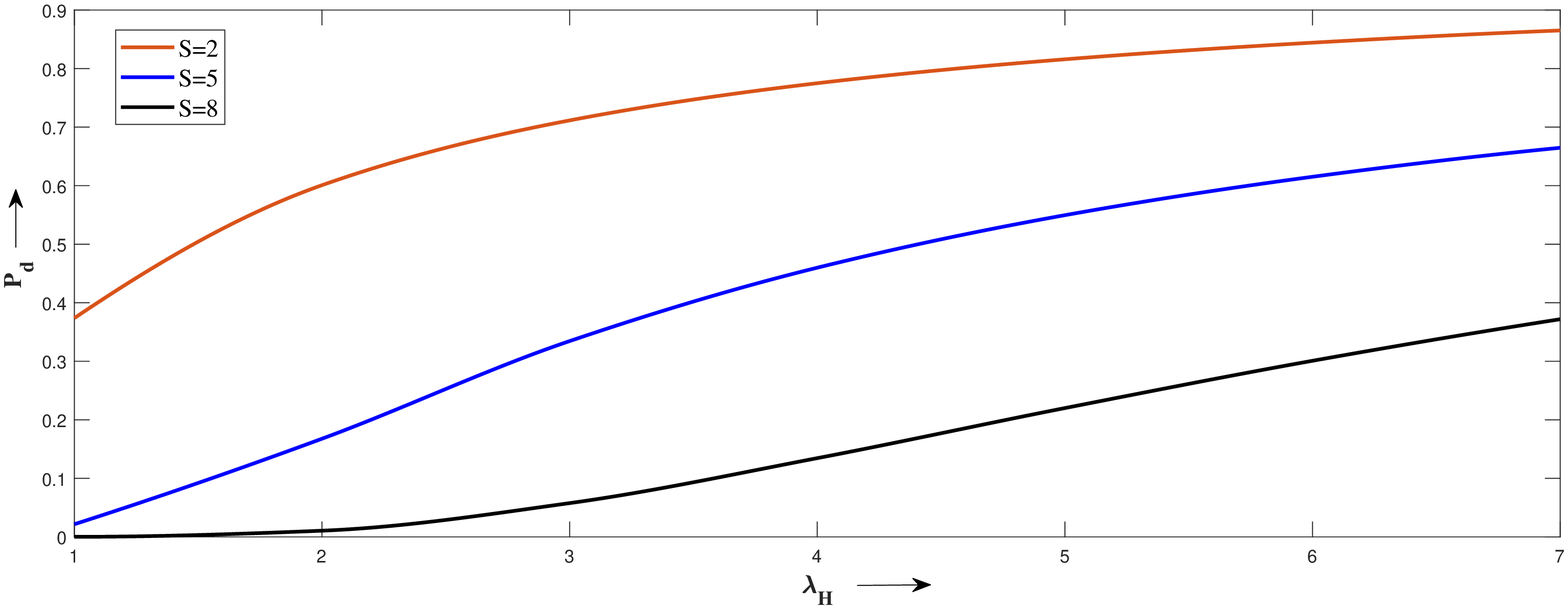}
		\label{fig:Pd1}}%
	\subfigure[$P_{preempt}$ versus  $\lambda_{\mathcal{H}}$]
	{\includegraphics[trim= 2.5cm 0.05cm 2cm 0.3cm, clip=true, height = 5.5cm,width = 0.5\textwidth]{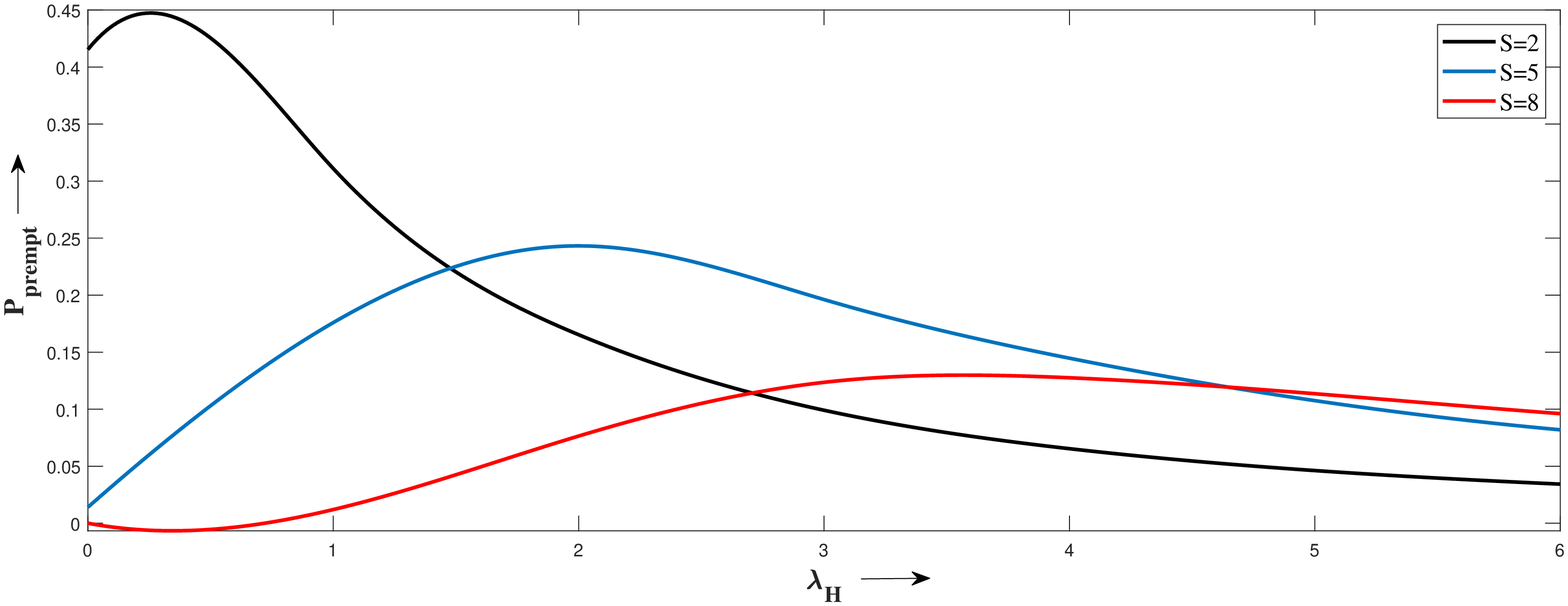}
		\label{fig:Pb1}}%
	\caption{Dependence of the dropping probability $P_d$ and preemption probability $P_{preempt}$ over arrival rate of a handoff call $\lambda_{\mathcal{H}}$. }
	\label{fig:Pb_Pd}
\end{figure}

Figures \ref{fig:Pd1} and \ref{fig:Pb1} represent $P_d$ and $P_{preempt}$ as functions of total number of channels $S$ and  $\lambda_{\mathcal{H}}$.  It can be observed from Figure \ref{fig:Pd1} that $P_d$ increases with respect to $\lambda_{\mathcal{H}}$ for a fixed value of $S$. Moreover, under the same value of $\lambda_{\mathcal{H}}$, $P_d$ decreases with respect to $S$. When  handoff calls arrive frequently in the system, all the channels are most likely to be occupied by handoff calls and consequently extra handoff calls will be dropped.
Therefore,  the value of $P_d$ increases with increasing $\lambda_{\mathcal{H}}$. In this scenario, if the total channels are increased in the system, more handoff call will be able to obtain the service, as a result, $P_d$ decreases.
Figure \ref{fig:Pb1} exhibits the impact of  $\lambda_{\mathcal{H}}$ and $S$ over $P_{preempt}$.
It is seen from the graph that the values of $P_{preempt}$ for different $S$, first increase, and then decrease. The cause for this behavior of $P_{preempt}$ lies in the following explanation.
When $\lambda_{\mathcal{H}}$  is relatively small, an arriving handoff call often finds at least one channel available, and consequently the ongoing service of a new call is not preempted by the arriving handoff call. As $\lambda_{\mathcal{H}}$ increases, the number of handoff calls also increase in the system. If an  arriving handoff call finds all the channels occupied and at least one of them is serving a new call, the service of that new call will be preempted by the arriving handoff call. Hence, $P_{preempt}$ increases and reaches maximum at some value of $\lambda_{\mathcal{H}}$. Further, the decreasing behaviour of $P_{preempt}$  is explained by the fact that, with the increment in $\lambda_{\mathcal{H}}$,  all the channels are occupied with handoff calls. Thus, the number of new calls in the service decreases and the probability that an arriving handoff call preempts the service of a new call decreases. 
These figures provide some advantageous results for the proposed model and  these results are opted to formulate an optimization problem  later on.

\begin{figure}
	\centering
	\subfigure[$P_d$ versus  $\lambda_{\mathcal{N}}$]
	{\includegraphics[trim= 3cm 0.1cm 3.8cm 0.4cm, height = 5.55cm,width = 0.5\textwidth]{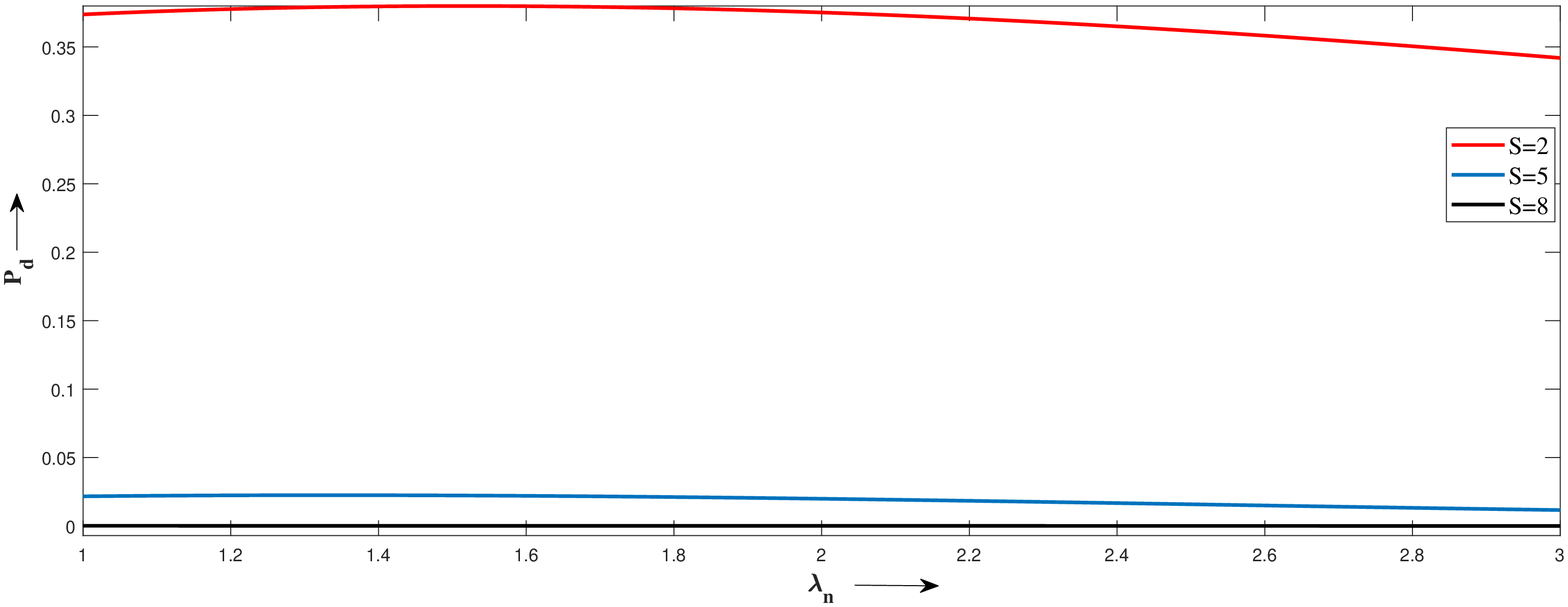}
		\label{fig:Pdn}}%
	\subfigure[$P_{preempt}$ versus  $\lambda_{\mathcal{N}}$]
	{\includegraphics[trim= 2.5cm 0.05cm 2cm 0.3cm, clip=true, height = 5.5cm,width = 0.5\textwidth]{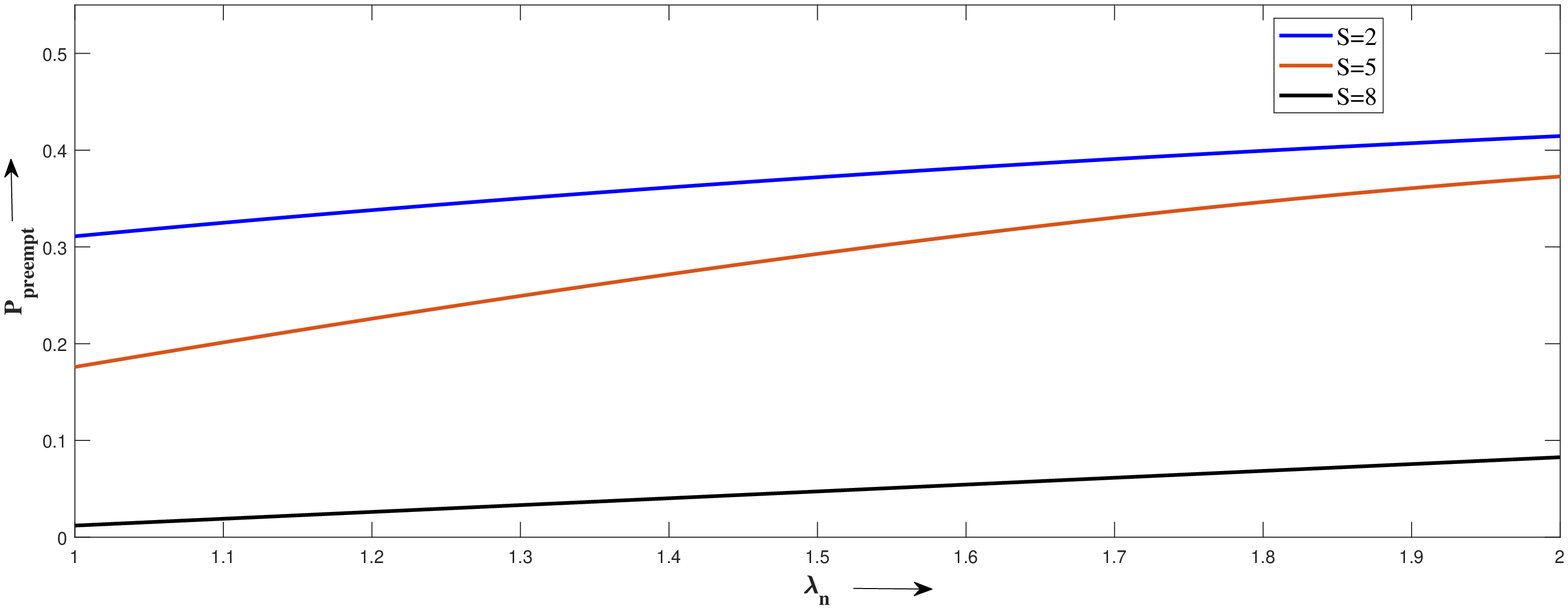}
		\label{fig:Pbn}}%
	\caption{Dependence of the dropping probability $P_d$ and preemption probability $P_{preempt}$ over arrival rate of a new call $\lambda_{\mathcal{N}}$. }
	\label{fig:Pb_Pd}
\end{figure}

Figures \ref{fig:Pdn} and \ref{fig:Pbn} describe the behaviour of $P_d$ and $P_{preempt}$ corresponding to  $\lambda_{\mathcal{N}}$ and $S$. The  plot \ref{fig:Pdn} exhibits that $P_d$ is less  affected by an increment in $\lambda_{\mathcal{N}}$ but $P_d$ decreases  with the increasing value of $S.$
Whereas, an obvious increasing  behaviour  of $P_{preempt}$ is observed from Figure \ref{fig:Pbn} when $\lambda_{\mathcal{N}}$  increases. More interestingly, $P_{preempt}$ decreases when $S$  increases in the system. Such  behaviour of $P_{preempt}$ can be explained as follows. When the number of new calls receiving  service in the system is very less, the probability of preemption for new calls also decreases and  as $\lambda_{\mathcal{N}}$ increases, $P_{preempt}$ also increases. Though, this probability can be reduced by increasing $S$ in the system.\\

\textbf{Experiment 2:} In this experiment, the impact of the service rate of handoff calls ($\mu_{\mathcal{H}}$) is shown over $P_d$ and $P_{preempt}$.

\begin{figure}
	\centering
	\subfigure[$P_d$ versus  $\mu_{\mathcal{H}}$]
	{\includegraphics[trim= 3cm 0.1cm 3.8cm 0.4cm, height = 5.55cm,width = 0.5\textwidth]{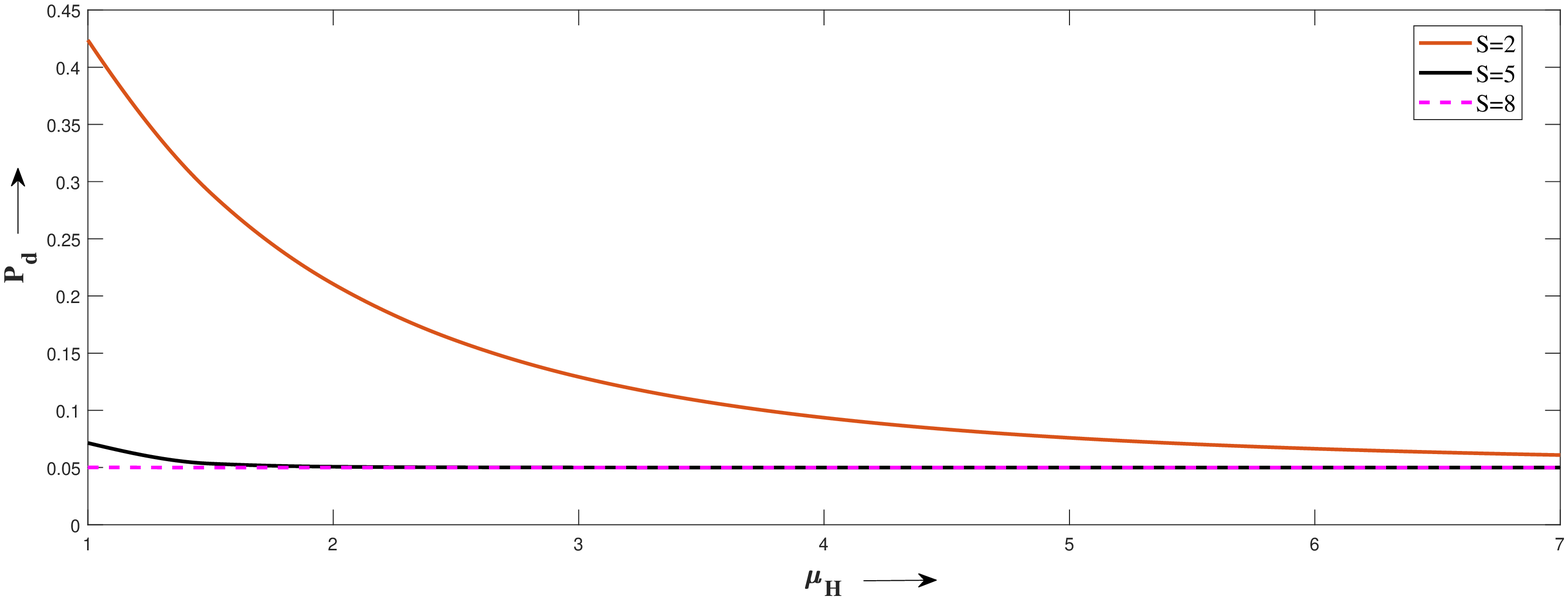}
		\label{fig:Pd2}}%
	\subfigure[$P_{preempt}$ versus  $\mu_{\mathcal{H}}$]
	{\includegraphics[trim= 2.5cm 0.05cm 2cm 0.3cm, clip=true, height = 5.5cm,width = 0.5\textwidth]{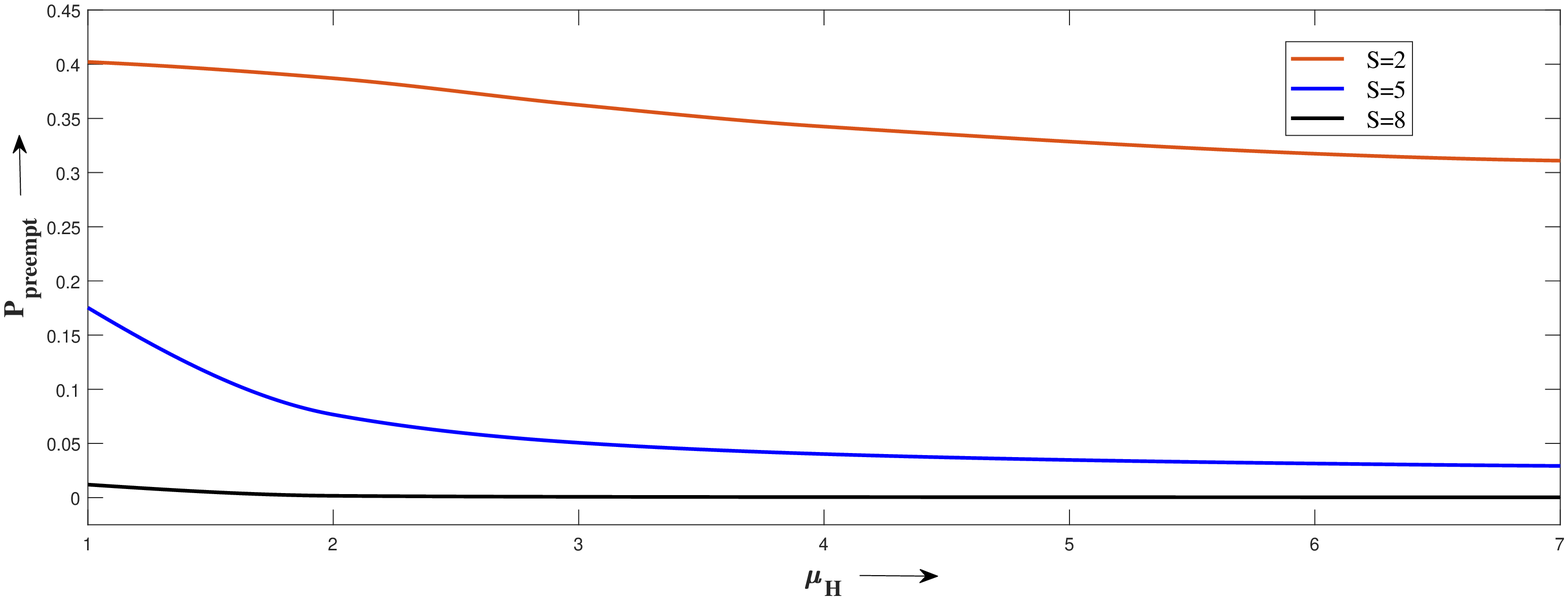}
		\label{fig:Pb2}}%
	\caption{Dependence of the dropping probability $P_d$ and preemption probability $P_{preempt}$ over service rate of a handoff call $\mu_{\mathcal{H}}$. }
	\label{fig:Pb_Pd}
\end{figure}

Figures \ref{fig:Pd2} and \ref{fig:Pb2} show the decreasing behaviour of  $P_d$ and $P_{preempt}$ as  $\mu_{\mathcal{H}}$ increases in the system. Further, it can also be observed from  Figures \ref{fig:Pd2} and \ref{fig:Pb2} that, for the fixed value of $\mu_{\mathcal{H}}$,  $P_d$ and $P_{preempt}$  decrease as  $S$ increases. An intuitive explanation for this finding can easily be given as follows. As $\mu_{\mathcal{H}}$ increases, the calls are served with increasing rate, consequently the loss probabilities for both types of calls  decrease. If $S$ is increased in the system, the chances  for both types of calls to obtain the service also increase and consequently $P_d$ and $P_{preempt}$ decrease.
\\

\textbf{Experiment 3:} The main purpose of this experiment is to observe the behaviour of intensity by which a retrial call is successfully connected to an available channel ($\theta_r^{succ}$) and the probability that a retrial call will exit the system without obtaining the service ($P_{leave}^{no-service}$) with respect to the retrial rate ($\theta$).

\begin{figure}
	\centering
	\subfigure[$\theta_r^{succ}$ versus  $\theta$]
	{\includegraphics[trim= 2cm 0.1cm 3cm 0.4cm, height = 5.5cm,width = 0.5\textwidth]{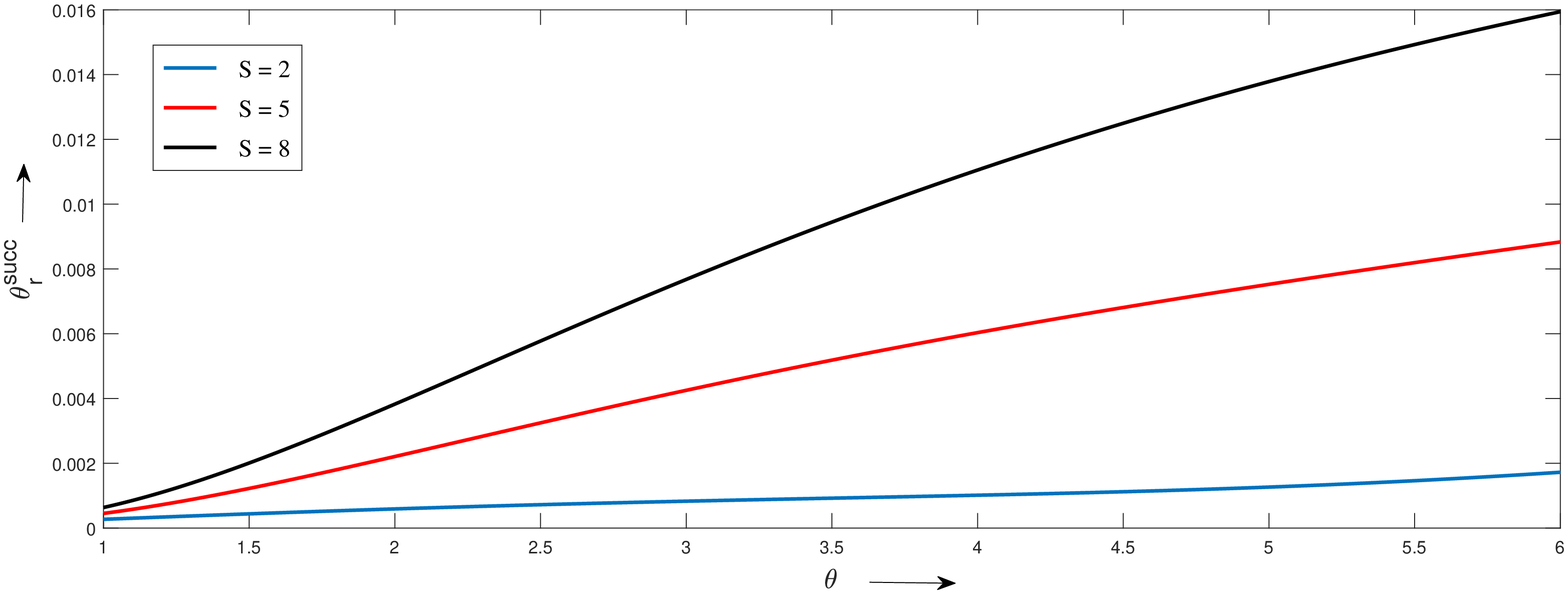}
		\label{fig:rate_capture_channel}}%
	\subfigure[$P_{leave}^{no-service}$ versus $\theta$]
	{\includegraphics[trim= 2.5cm 0.05cm 2cm 0.5cm, clip=true, height = 5.5cm,width = 0.5\textwidth]{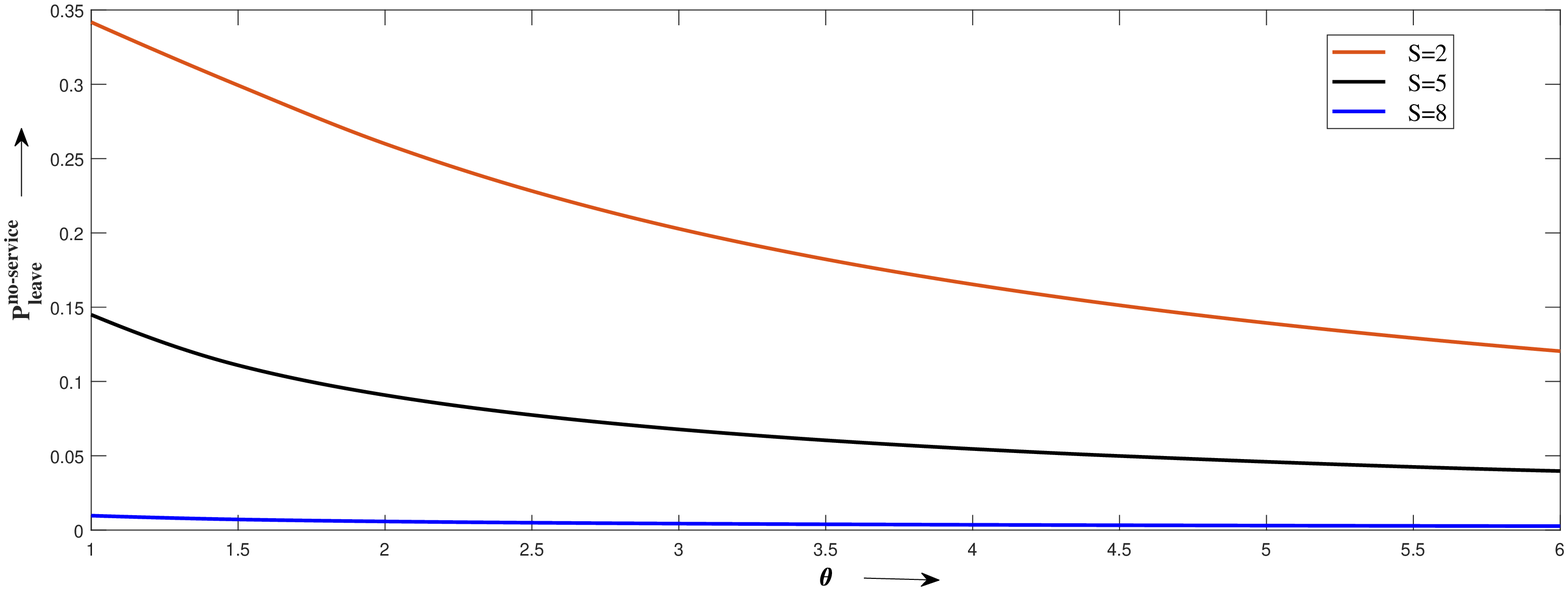}
		\label{fig:impatience_retrial}}%
	\caption{Dependence of the intensity by which a retrial call is successfully connected to an available channel $\theta_r^{succ}$ and    the probability that a retrial call will exit the system without obtaining the service $P_{leave}^{no-service}$ over retrial rate $\theta$. }
	\label{fig:retrial_rate}
\end{figure}

Figures \ref{fig:rate_capture_channel} and \ref{fig:impatience_retrial} represent the behaviour of $\theta_r^{succ}$ and $P_{leave}^{no-service}$ with respect to  $\theta$, respectively.  It can be seen from the graphs when the value of $\theta$ increases,  $\theta_r^{succ}$ increases. This impact can easily be explained as follows. When  $\theta$ increases,  the probability of a retrial call getting a connection  also increases. For the fixed value of $\theta$, this probability is substantially greater for the large value of $S.$ The opposite impact of  $\theta$ over $P_{leave}^{no-service}$ can be observed from Figure \ref{fig:impatience_retrial}. When $\theta$ increases, the probability that a retrial call  leaves the system without obtaining the service  decreases. Moreover, decreasing behaviour of  $P_{leave}^{no-service}$ is observed with respect to $S$ for the fixed value of $\theta$. When the number of  available channels increases in the system, the probability for retrial call to obtain the service also increases. From  \ref{fig:Pd1} and \ref{fig:Pb1}, it can be observed that $\lambda_{\mathcal{H}}$ has vital impact over $P_d$ and $P_{preempt}$ when all other parameters of the system are fixed according to the system requirement. This observation is the main motivation for the formulation of the  optimization problem illustrated in Section \ref{section6}. 
%##############################

%####################################
%##############################

%####################################
%##############################

%########################
%######################

%#################################### 

\section{Optimization Problem for Traffic Control} \label{section6}

\noindent The loss probabilities are considered as performance determining factors for  cellular networks. Therefore,  $P_d$ and $P_{preempt}$ are considered crucial factors in the proposed model. An increment in $P_d$  as well as in $P_{preempt}$ indicates unsatisfactory level of service for the customers. Thus, it is required to find the optimal values of parameters in such a way that the loss probabilities should not exceed some pre-defined values. In Section \ref{section5}, a detailed analysis of loss probabilities with respect to several parameters has been provided. It can be observed from the results that $P_d$ and $P_{preempt}$ are mostly affected by the arrival rate of handoff calls $\lambda_{\mathcal{H}}$ and the total number of channels $S.$ The service provider certainly cannot determine $\lambda_{\mathcal{H}}$, yet an approximated value of $\lambda_{\mathcal{H}}$ can be estimated  in order to keep sufficient channels to provide service. Thus, a non-trivial optimization problem is proposed with the decision variables $S$ and $ \lambda_{\mathcal{H}}$ given as:\\
\begin{center}
$\begin{array}{lll}
&\textrm{min } &S\\
&\textrm{subject to},& P_d( S,\lambda_{\mathcal{H}}) \leq \epsilon_1,\\
&&P_{preempt}(S,\lambda_{\mathcal{H}}) \leq \epsilon_2, \\
&& S,\lambda_{\mathcal{H}} \geq 0.
\end{array}$\\
\end{center}
Here, $\epsilon_1$ and $\epsilon_2$ are pre-defined values depending on the tolerance of the system for dropping probability and preemption probability, respectively.
The above defined  constraints are non-linear and highly complex.  Thus heuristic approaches have been used to obtain the optimal solution such as  Direct Search method ($DS$), Particle Swarm Optimization method (\textit{PSO}) and Simulated Annealing method ($SA$). Let the objective function be denoted as $F$  and the optimal solution as $\hat{x} = \{S^*,\lambda_{\mathcal{H}}^*\}$ for the  algorithms. Consider $\epsilon_1=10^{-4}$ and $\epsilon_2=10^{-4}$  for the numerical computation of the proposed optimization problem. The detailed sensitivity analysis of the optimization problem is provided as follows.

\subsection{Direct Search Method}

The \textit{DS} method can be applied to non-linear, convex or non-convex optimization problems. It does not need to compute the gradient of the objective function like gradient based  optimization methods (\cite{torczon1998evolutionary}).
Due to the easy implementation and computation, this method is used widely to obtain the optimal solution for any type of optimization problem. The algorithm for this method works as follows:\\

\textbf{Algorithm 1:}
\begin{itemize}
	\item Step 1: Fix $\lambda_{\mathcal{N}}$, $\mu_{\mathcal{H}}$, $\mu_{\mathcal{N}}$, $\theta$, $\epsilon_1$ and $\epsilon_2.$ Initialize $S=2.$
	\item Step 2: Find the values of $\lambda_{\mathcal{H}}^1$ and $\lambda_{\mathcal{H}}^2$ which satisfy $P_d( S,\lambda_{\mathcal{H}}^1) \leq \epsilon_1$ and $P_{preempt}(S,\lambda_{\mathcal{H}}^2) \leq \epsilon_2$, respectively.
	\item Step 3: If $\lambda_{\mathcal{H}}^1$ = $\lambda_{\mathcal{H}}^2$ = $\lambda_{\mathcal{H}}^*$, declare $S$ and $\lambda_{\mathcal{H}}^*$ as optimal values. Otherwise put $S = S+1$ and repeat Step 2.
\end{itemize}
Table \ref{tab:my_label} represents the optimal values of $S^*$ and $\lambda_{\mathcal{H}}^*$  obtained by applying \textit{DS} method for different values of $\mu_{\mathcal{H}}$ and $\lambda_{\mathcal{N}}$.    
\begin{table}[]
	\centering
	\scalebox{0.7}
	{
		\begin{tabular}{llllllll}
			\hline   
			$\lambda_{\mathcal{N}} = 0.1$,    $\mu_{\mathcal{N}}=1$ & & & & & & &\\
			\hline        
			$\mu_{\mathcal{H}}$ & 0.5 & 0.625 & 0.75 & 0.875 & 1 & 1.125 & 1.25 \\
			$S^*$ & 3 &3 &3 & $-$& $-$& 4& 4\\
			$\lambda_{\mathcal{H}}^*$ & 0.2 & 0.25 & 0.3 & $-$& $-$&0.925 & 1.05\\
			$P_d^*$ & $8.27 \times 10^{-4}$ & $7.75 \times 10^{-4}$ & $7.40\times 10^{-4}$ & $-$& $-$& $9.28 \times 10^{-4}$ & $9.19 \times 10^{-4}$
			\\
			$P_{preempt}^*$ & $ 7.65 \times 10^{-4}$ & $8.33 \times 10^{-4}$ & $9.03 \times 10^{-4}$ & $-$ & $-$ & $7.87 \times 10^{-4}$ & $8.16 \times 10^{-4}$
			\\
			\hline
			$\lambda_{\mathcal{N}} = 0.2$,    $\mu_{\mathcal{N}}=1$ & & & & & & &\\
			\hline        
			$\mu_{\mathcal{H}}$ & 0.5 & 0.625 & 0.75 & 0.875 & 1 & 1.125 & 1.25 \\
			$S^*$ & 4 & $-$ &4 & 4& 4& 4& 4\\
			$\lambda_{\mathcal{H}}^*$ & 0.25 & $-$ & 0.575 & 0.675& 0.425 & 0.875 & 0.975\\
			$P_d^*$ & $6.40 \times 10^{-4}$ & $-$ & $8.85\times 10^{-4}$ & $8.81\times 10^{-4}$& $8.28\times 10^{-4}$& $8.75 \times 10^{-4}$ & $9.43 \times 10^{-4}$
			\\
			$P_{preempt}^*$ & $ 7.81 \times 10^{-4}$ & $-$ & $7.35 \times 10^{-4}$ & $7.29\times 10^{-4}$ & $9.04\times 10^{-4}$ & $7.40 \times 10^{-4}$ & $8.41 \times 10^{-4}$
			\\
			\hline
			$\lambda_{\mathcal{N}} = 0.3$,    $\mu_{\mathcal{N}}=1$ & & & & & & &\\
			\hline        
			$\mu_{\mathcal{H}}$ & 0.5 & 0.625 & 0.75 & 0.875 & 1 & 1.125 & 1.25 \\
			$S^*$ & 4 & 4 & 4 & $-$ & 4 & 4& 4\\
			$\lambda_{\mathcal{H}}^*$ & 0.35 & 0.45 & 0.575 & $-$ & 0.775 &0.875 & 0.975\\
			$P_d^*$ & $7.52 \times 10^{-4}$ & $7.81 \times 10^{-4}$ & $9.95\times 10^{-4}$ & $-$ & $9.63\times 10^{-4}$ & $9.52 \times 10^{-4}$ & $9.42 \times 10^{-4}$
			\\
			$P_{preempt}^*$ & $ 5.90 \times 10^{-4}$ & $6.49 \times 10^{-4}$ & $8.30 \times 10^{-4}$ & $
			-$ & $9.11\times 10^{-4}$ & $9.53 \times 10^{-4}$ & $9.94 \times 10^{-4}$
			\\
			\hline
			$\lambda_{\mathcal{N}} = 0.4$,    $\mu_{\mathcal{N}}=1$ & & & & & & &\\
			\hline        
			$\mu_{\mathcal{H}}$ & 0.5 & 0.625 & 0.75 & 0.875 & 1 & 1.125 & 1.25 \\
			$S^*$ & 4 &$-$ &4 & $-$& 4& $-$& 5\\
			$\lambda_{\mathcal{H}}^*$ & 0.35 & $-$ & 0.525 & $-$ & 0.7 & $-$ & 1.475\\
			$P_d^*$ & $8.70 \times 10^{-4}$ & $-$ & $7.12\times 10^{-4}$ & $-$& $6.32\times 10^{-4}$ & $-$ & $8.08 \times 10^{-4}$
			\\
			$P_{preempt}^*$ & $ 9.80 \times 10^{-4}$ & $-$ & $9.59 \times 10^{-4}$ & $-$ & $9.97 \times 10^{-4}$ & $-$ & $9.29 \times 10^{-4}$
			\\
			\hline
			$\lambda_{\mathcal{N}} = 0.5$,    $\mu_{\mathcal{N}}=1$ & & & & & & &\\
			\hline        
			$\mu_{\mathcal{H}}$ & 0.5 & 0.625 & 0.75 & 0.875 & 1 & 1.125 & 1.25 \\
			$S^*$ & 5 & $-$ & 6 & 6 & 6 & 6 & $-$ \\
			$\lambda_{\mathcal{H}}^*$ & 0.55 & $-$ & 1.25 & 1.45 & 1.475 & 1.9 & $-$\\
			$P_d^*$ & $7.94 \times 10^{-4}$ & $-$ & $9.47\times 10^{-4}$ & $8.76 \times 10^{-4}$ &$9.79\times 10^{-4}$& $9.38 \times 10^{-4}$ & $-$
			\\
			$P_{preempt}^*$ & $ 9.21 \times 10^{-4}$ & $-$ & $9.29 \times 10^{-4}$ & $9.14\times 10^{-4}$ & $9.96\times 10^{-4}$ & $9.77 \times 10^{-4}$ & $-$
			\\
			\hline
			$\lambda_{\mathcal{N}} = 0.6$,   $\mu_{\mathcal{N}}=1$ & & & & & & &\\
			\hline        
			$\mu_{\mathcal{H}}$ & 0.5 & 0.625 & 0.75 & 0.875 & 1 & 1.125 & 1.25 \\
			$S^*$ & 6 & $-$ & 6 & 6 & 7 & $-$ & 7\\
			$\lambda_{\mathcal{H}}^*$ & 0.775 & $-$ & 1.175 & 1.35 & 1.6 & $-$ & 2.1\\
			$P_d^*$ & $7.36 \times 10^{-4}$ & $-$ & $6.70\times 10^{-4}$ & $7.31 \times 10^{-4}$ & $7.03 \times 10^{-4}$ & $-$ & $7.33 \times 10^{-4}$
			\\
			$P_{preempt}^*$ & $ 8.91 \times 10^{-4}$ & $-$ & $9.11 \times 10^{-4}$ & $8.50\times 10^{-4}$ & $8.84\times 10^{-4}$ & $-$ & $8.48 \times 10^{-4}$
			\\
			\hline
	\end{tabular}}
	\caption{Optimal values of $S^*$ and $\lambda_{\mathcal{H}}^*$ for different values of $\mu_{\mathcal{H}}$ and $\lambda_{\mathcal{N}}$ by applying DS Method}
	\label{tab:my_label}
\end{table}

\subsection{Particle Swarm Optimization Method}

The algorithm for \textit{PSO} method was developed by \cite{eberhart1995particle}. In this method, a population of particles are explored and exploited by applying search techniques in different dimensions.  Each particle has velocity and position which defines a particle's best known local and global positions in the feasible space. 
The \textit{PSO} method is suitable for non-differentiable constrained and unconstrained optimization problems as it does not require the computation of the gradient. This technique can search discrete and continuous decision variables at the same time. Since, the proposed optimization problem is constrained optimization problem, penalty function approach has been applied to convert the same into an unconstrained optimization problem.
The algorithm for \textit{PSO} method works as follows:\\

\textbf{Algorithm 2:}
\begin{itemize}
	\item Step 1: Set parameters of \textit{PSO}, i.e., $w_{min}$, $w_{max}$, $c_{1}$ and $c_{2}$; where $w \in [w_{min},w_{max}]$ is an inertia factor and $c_1,c_2$ are acceleration factors. Fix the parameters,  $\lambda_{\mathcal{N}}$, $\mu_{\mathcal{H}}$, $\mu_{\mathcal{N}}$, $\theta$, $\epsilon_1$ and $\epsilon_2$ for the computation of the proposed optimization problem as per the system requirement.
	\item Step 2: Initialize population of particles having positions $x$, velocities $V$ and maximum iterations $maxite$.
	\item Step 3: Convert a constrained optimization problem to an unconstrained optimization problem by using penalty function approach. Compute fitness function of particles $F_i = F(x_i),\forall i$; where $i$ is population of particles and $x = \{x^{(1)},x^{(2)}\} = \{S,\lambda_{\mathcal{H}}\}$.
	\item Step 4: Initialize the partial best solution $P_B = x$ and the global best solution $G_B = arg~min\{F(x);x \in P_B\}$.
	\item Step 5: Set $Iteration=1$. Generate two random numbers $U_1 = U(0,1)$ and $U_2 = U(0,1)$.
	\item Step 6: Update the particle's positions and velocity as   $x = x+V$ and $V = wV + c_1 U_1 (P_B-x)+ c_2 U_2(G_B-x)$.
	\item Step 7: Update the partial best solution $P_B$ and $G_B$.
	\item Step 8: Repeat Steps 5-7 until $Iteration < maxite$.
	\item Step 9: Obtain the optimal solution $\hat{x} = G_B$ and $F^* = F(\hat{x})$.
\end{itemize}
For various values of $\mu_{\mathcal{H}}$ and $\lambda_{\mathcal{N}}$, the \textit{PSO} method is executed with initial 60 random generated particles, $w_{max}=0.9$, $w_{min}=0.4$, $c_1=c_2=2$, and $maxite=200.$ 
Table \ref{tab:my_label1} exhibits the optimal values of $S^*$ and $\lambda_{\mathcal{H}}^*$  obtained by applying \textit{PSO} for different combinations of $\mu_{\mathcal{H}}$ and $\lambda_{\mathcal{N}}$. 

\subsection{Simulated Annealing  Method}
$SA$ is a meta-heuristic approach which has the convergence rate  comparatively slower than other heuristic methods, but this method does not stuck in the pool of local optimal solutions. Thus, this method guarantees to  obtain the global optimal solution successfully \cite{kirkpatrick1983optimization}. This method is applicable for complex non-differentiable constrained and unconstrained optimization problems. The algorithm for $SA$ method works as follows:\\

\textbf{Algorithm 3:}
\begin{itemize}
	\item Step 1: Fix the parameters $\lambda_{\mathcal{N}}$, $\mu_{\mathcal{H}}$, $\mu_{\mathcal{N}}$, $\theta$, $\epsilon_1$ and $\epsilon_2$  for the computation of the objective function. Convert a constrained optimization problem to an unconstrained optimization problem by using penalty function approach.
	
	\item Step 2: Initialize $x$ as the initial state. Generate $x^{'}$ $= x + \Delta x$ a neighbor state of $x,$ where  $x = \{x^{(1)},x^{(2)}\} = \{S,\lambda_{\mathcal{H}}\}$.
	\item Step 3: Compute $\Delta F(x) = F(x) - F(x^{'}).$
	\item Step 4: Generate $P = e^{(-\Delta F(x)/T)}$, the acceptance probability of $x,$ where $T$, a temperature parameter,  is  evaluated randomly by computing the mean of different values of given objective function.
	\item Step 5: Generate $R$ = $U(0,1)$, the acceptance probability of $x^{'}$.
	\item Step 6: If $\Delta F(x) < 0$; $x$ is the actual state else if $R>P$; $x^{'}$ is the actual state.
	\item Step 7: Repeat Steps 2-6 until $\Delta F(x) < 10^{-6}$.
	\item Step 8: Obtain the optimal solution  $\hat{x}  =\{S^*,\lambda_{\mathcal{H}}^*\}$.

\end{itemize}

Table \ref{tab:my_label2} exhibits the optimal values of $S^*$ and $\lambda_{\mathcal{H}}^*$  obtained by employing \textit{SA} method for different combinations of $\mu_{\mathcal{H}}$ and $\lambda_{\mathcal{N}}$. 

All the results were obtained by MATLAB software, which were run on a computer with Intel Core i7-6700 3.40GHz CPU and 16 GB of RAM. It can be observed from Table \ref{tab:my_label}, Table \ref{tab:my_label1} and Table \ref{tab:my_label2} that \textit{PSO} method and \textit{SA} method are more efficient. These methods provide optimal solution for each possible combination of $\mu_{\mathcal{H}}$ and $\lambda_{\mathcal{N}}$. Though, the number of iterations in \textit{SA} method are more than the number of iterations in \textit{PSO} method, \textit{SA} method always provide the global optimal solution.  For this optimization problem both the methods are providing identical results which determines the reliability of  the optimal solutions.

%The obtained results represent the efficiency of \textit{PSO} method and \textit{SA} method over \textit{DS} method.  It can be observed from the Table \ref{tab:my_label} that \textit{DS} method fails to provide optimal solutions for some values of $\lambda_{\mathcal{N}}$ and $\mu_{\mathcal{H}}$. While,  \textit{PSO} and \textit{SA} methods provide optimal solutions for each combination of $\mu_{\mathcal{H}}$ and $\lambda_{\mathcal{N}}$. In  \textit{SA} method, the number of iterations to provide optimal solutions are more in comparison to \textit{PSO} method but \textit{SA} method  guarantees for the global optimal solution. Whereas, \textit{PSO} method is inclined to trap into the pool of local optimal solutions and consequently does not guarantee for the global optimal solution. The application of \textit{SA} method ensures that the result obtained by  \textit{PSO} method are the global ones for the proposed optimization problem. These results are almost identical which determines the reliability of the optimal solutions. %Therefore, the obtained results demonstrate the efficiency of  \textit{PSO} and \textit{SA} methods.

\begin{table}[]
	\centering
	\scalebox{0.7}
	{
	\begin{tabular}{llllllll}
		\hline   
		$\lambda_{\mathcal{N}} = 0.1$,    $\mu_{\mathcal{N}}=1$ & & & & & & &\\
		\hline
		Iterations & 101 & 103 &102 & 103 & 103 & 103 & 101\\
		$\mu_{\mathcal{H}}$ & 0.5 & 0.625 & 0.75 & 0.875 & 1 & 1.125 & 1.25 \\
		$S^*$ & 2 &2 &2 &2& 2& 2& 2\\
		$\lambda_{\mathcal{H}}^*$ & 0.013 & 0.01 & 0.01 & 0.0324 & 0.01 & 0.01 & 0.01\\
		$P_d^*$ & $3.28 \times 10^{-6}$ & $7.699 \times 10^{-7}$ & $3.887 \times 10^{-7}$ & $8.921 \times 10^{-5}$& $1.3138 \times 10^{-7}$ & $8.382 \times 10^{-8}$ & $5.5967 \times 10^{-8}$
		\\
		$P_{preempt}^*$ & $2.87 \times 10^{-5}$ & $1.974 \times 10^{-5}$ & $1.767 \times 10^{-5}$ & $1.6327 \times 10^{-4}$ & $1.5695 \times 10^{-5}$ & $1.517 \times 10^{-5}$ & $1.4809 \times 10^{-5}$
		\\
		\hline
		$\lambda_{\mathcal{N}} = 0.2$,    $\mu_{\mathcal{N}}=1$ & & & & & & &\\
		\hline
		Iterations & 104 & 110 & 103 & 103 & 103 & 103 & 103\\
		$\mu_{\mathcal{H}}$ & 0.5 & 0.625 & 0.75 & 0.875 & 1 & 1.125 & 1.25 \\
		$S^*$ & 2 &2 & 2 &2& 2& 2& 2\\
		$\lambda_{\mathcal{H}}^*$ & 0.01 & 0.01 & 0.01001 & 0.01001 & 0.01001 & 0.01001 & 0.01001\\
		$P_d^*$ & $1.14 \times 10^{-6}$ & $2.17 \times 10^{-6}$ & $1.14 \times 10^{-7}$ & $6.571 \times 10^{-7}$& $4.042 \times 10^{-7}$ & $2.6624\times 10^{-7}$ & $3.27 \times 10^{-7}$
		\\
		$P_{preempt}^*$ & $1.87 \times 10^{-4}$ & $1.98 \times 10^{-4}$ & $7.37\times 10^{-7}$ & $1.8023 \times 10^{-4}$ & $1.75 \times 10^{-4}$ & $1.7271 \times 10^{-4}$ & $4.75 \times 10^{-4}$
		\\
		\hline
		$\lambda_{\mathcal{N}} = 0.3$,    $\mu_{\mathcal{N}}=1$ & & & & & & &\\
		\hline
		Iterations & 102 & 102 & 103 & 103 & 103 & 103 & 103\\
		$\mu_{\mathcal{H}}$ & 0.5 & 0.625 & 0.75 & 0.875 & 1 & 1.125 & 1.25 \\
		$S^*$ & 2 & 2 & 2 & 2 & 2 &  2 & 2\\
		$\lambda_{\mathcal{H}}^*$ & 0.01 & 0.01 & 0.01 & 0.01 & 0.01 & 0.01 & 0.01\\
		$P_d^*$ & $7.89 \times 10^{-6}$ & $3.79 \times 10^{-6}$ & $2.047 \times 10^{-6}$ & $1.20 \times 10^{-6}$& $7.53 \times 10^{-7}$ & $4.95 \times 10^{-7}$ & $3.39 \times 10^{-7}$
		\\
		$P_{preempt}^*$ & $8.15 \times 10^{-4}$ & $7.65 \times 10^{-4}$ & $7.36 \times 10^{-4}$ & $7.18 \times 10^{-4}$ & $7.06 \times 10^{-4}$ & $6.98 \times 10^{-4}$ & $6.92 \times 10^{-4}$
		\\
		\hline
		$\lambda_{\mathcal{N}} = 0.4$,    $\mu_{\mathcal{N}}=1$ & & & & & & &\\
		\hline
		Iterations & 102 & 102 &103 & 103 & 103 & 103 & 103\\
		$\mu_{\mathcal{H}}$ & 0.5 & 0.625 & 0.75 & 0.875 & 1 & 1.125 & 1.25 \\
		$S^*$ & 3 &3 &3 & 3 & 3& 3& 3\\
		$\lambda_{\mathcal{H}}^*$ & 0.0105 & 0.0105 & 0.0105 & 0.0105 & 0.01001 & 0.01001 & 0.01001\\
		$P_d^*$ & $1.21 \times 10^{-8}$ & $4.02 \times 10^{-9}$ & $1.60\times 10^{-9}$ & $7.24 \times 10^{-10}$& $1.66 \times 10^{-10}$ & $5.001 \times 10^{-9}$ & $9.45 \times 10^{-11}$
		\\
		$P_{preempt}^*$ & $2.53 \times 10^{-5}$ & $2.31 \times 10^{-5}$ & $2.19 \times 10^{-5}$ & $2.62 \times 10^{-5}$ & $2.68 \times 10^{-5}$ & $2.04 \times 10^{-5}$ & $2.02 \times 10^{-5}$
		\\
		\hline
		$\lambda_{\mathcal{N}} = 0.5$,    $\mu_{\mathcal{N}}=1$ & & & & & & &\\
		\hline
		Iterations & 102 & 102 &102 & 103 & 103 & 103 & 103\\
		$\mu_{\mathcal{H}}$ & 0.5 & 0.625 & 0.75 & 0.875 & 1 & 1.125 & 1.25 \\
		$S^*$ & 3 &3 &3 &3& 3& 3& 3\\
		$\lambda_{\mathcal{H}}^*$ & 0.0105 & 0.0105 & 0.0105 & 0.0105 & 0.0105 & 0.0105 & 0.0105\\
		$P_d^*$ & $1.76 \times 10^{-8}$ & $6.11 \times 10^{-9}$ & $2.50 \times 10^{-9}$ & $1.15 \times 10^{-9}$& $5.82 \times 10^{-10}$ & $3.16 \times 10^{-10}$ & $2.69 \times 10^{-10}$
		\\
		$P_{preempt}^*$ & $7.40 \times 10^{-5}$ & $6.93 \times 10^{-5}$ & $9.24 \times 10^{-5}$ & $6.67 \times 10^{-5}$ & $6.39 \times 10^{-5}$ & $6.34 \times 10^{-5}$ & $8.73 \times 10^{-5}$
		\\
		\hline
		$\lambda_{\mathcal{N}} = 0.6$,    $\mu_{\mathcal{N}}=1$ & & & & & & &\\
		\hline
		Iterations & 102 & 102 &103 & 103 & 102 & 102 & 102\\
		$\mu_{\mathcal{H}}$ & 0.5 & 0.625 & 0.75 & 0.875 & 1 & 1.125 & 1.25 \\
		$S^*$ & 3 & 3 &3 & 3 & 3& 3 & 3\\
		$\lambda_{\mathcal{H}}^*$ & 0.0105 & 0.0105 & 0.0105 & 0.0105 & 0.0105 & 0.0105 & 0.0105\\
		$P_d^*$ & $2.35 \times 10^{-8}$ & $8.39 \times 10^{-9}$ & $3.51 \times 10^{-9}$ & $1.65 \times 10^{-9}$& $8.49 \times 10^{-10}$ & $4.67 \times 10^{-10}$ & $2.71 \times 10^{-10}$
		\\
		$P_{preempt}^*$ & $1.75 \times 10^{-4}$ & $1.66 \times 10^{-4}$ & $1.61 \times 10^{-4}$ & $1.58 \times 10^{-4}$ & $1.56 \times 10^{-4}$ & $1.55 \times 10^{-4}$ & $1.54 \times 10^{-4}$
		\\
		\hline
	\end{tabular}}
	\caption{Optimal values of $S^*$ and $\lambda_{\mathcal{H}}^*$ for different values of $\mu_{\mathcal{H}}$ and $\lambda_{\mathcal{N}}$ by applying PSO method}
	\label{tab:my_label1}
\end{table}
\begin{table}[H]
	\centering
	\scalebox{0.7}
	{
	\begin{tabular}{llllllll}
		\hline   
		$\lambda_{\mathcal{N}} = 0.1$,    $\mu_{\mathcal{N}}=1$ & & & & & & &\\
		\hline
		Iterations & 1000 & 1003 &1003 & 1003 & 1003 & 1003 & 1003\\
		$\mu_{\mathcal{H}}$ & 0.5 & 0.625 & 0.75 & 0.875 & 1 & 1.125 & 1.25 \\
		$S^*$ & 2 &2 &2 &2& 2& 2& 2\\
		$\lambda_{\mathcal{H}}^*$ & 0.01 & 0.01 & 0.01 & 0.01 & 0.01 & 0.01 & 0.01\\
		$P_d^*$ & $1.7414 \times 10^{-6}$
 & $5.435 \times 10^{-10}$ & $5.5704 \times 10^{-11}$ & $2.3087 \times 10^{-11}$ & $1.0715 \times 10^{-11}$ &  $5.4264 \times 10^{-12}$ & $2.9458 \times 10^{-12}$
		\\
		$P_{preempt}^*$ & $2.3650 \times 10^{-5}$ & $2.1316 \times 10^{-8}$ & $1.7443 \times 10^{-8}$ & $1.5404 \times 10^{-8}$ & $1.4219 \times 10^{-8}$ & $1.3468 \times 10^{-8}$ & $1.2960 \times 10^{-8}$
		\\
		\hline
		$\lambda_{\mathcal{N}} = 0.2$,    $\mu_{\mathcal{N}}=1$ & & & & & & &\\
		\hline
		Iterations & 1000 & 1003 &1003 & 1003 & 1003 & 1003 & 1003\\
		$\mu_{\mathcal{H}}$ & 0.5 & 0.625 & 0.75 & 0.875 & 1 & 1.125 & 1.25 \\
		$S^*$ & 2 &2 &2 &2& 2& 2& 2\\
		$\lambda_{\mathcal{H}}^*$ & 0.01001 & 0.01001 & 0.01001 & 0.01001 & 0.01001 & 0.01001 & 0.01001\\
		$P_d^*$ & $1.15 \times 10^{-6}$ & $2.18 \times 10^{-6}$ & $1.15 \times 10^{-7}$ & $6.6 \times 10^{-7}$& $4.07 \times 10^{-7}$ & $2.67\times 10^{-7}$ & $3.28 \times 10^{-7}$
		\\
		$P_{preempt}^*$ & $1.88 \times 10^{-4}$ & $1.99 \times 10^{-4}$ & $7.38\times 10^{-7}$ & $1.81 \times 10^{-4}$ & $1.76 \times 10^{-4}$ & $1.73 \times 10^{-4}$ & $4.76 \times 10^{-4}$
		\\
		\hline
		$\lambda_{\mathcal{N}} = 0.3$,    $\mu_{\mathcal{N}}=1$ & & & & & & &\\
		\hline
		Iterations & 1000 & 1003 &1003 & 1003 & 1003 & 1003 & 1003\\
		$\mu_{\mathcal{H}}$ & 0.5 & 0.625 & 0.75 & 0.875 & 1 & 1.125 & 1.25 \\
		$S^*$ & 2 & 2 & 2 & 2 & 2 &  2 & 2\\
		$\lambda_{\mathcal{H}}^*$ & 0.01001 & 0.01001 & 0.01001 & 0.01001 & 0.01001 & 0.01001 & 0.01001\\
		$P_d^*$ & $7.89 \times 10^{-6}$ & $3.79 \times 10^{-6}$ & $2.047 \times 10^{-6}$ & $1.20 \times 10^{-6}$& $7.53 \times 10^{-7}$ & $4.95 \times 10^{-7}$ & $3.39 \times 10^{-7}$
		\\
		$P_{preempt}^*$ & $8.15 \times 10^{-4}$ & $7.65 \times 10^{-4}$ & $7.36 \times 10^{-4}$ & $7.18 \times 10^{-4}$ & $7.06 \times 10^{-4}$ & $6.98 \times 10^{-4}$ & $6.92 \times 10^{-4}$
		\\
		\hline
		$\lambda_{\mathcal{N}} = 0.4$,    $\mu_{\mathcal{N}}=1$ & & & & & & &\\
		\hline
			Iterations & 1000 & 1003 &1003 & 1003 & 1003 & 1003 & 1003\\
		$\mu_{\mathcal{H}}$ & 0.5 & 0.625 & 0.75 & 0.875 & 1 & 1.125 & 1.25 \\
		$S^*$ & 3 &3 &3 & 3 & 3& 3& 3\\
		$\lambda_{\mathcal{H}}^*$ & 0.01001 & 0.01001 & 0.01001 & 0.01001 & 0.01001 & 0.01001 & 0.01001\\
		$P_d^*$ & $1.22 \times 10^{-8}$ & $4.03 \times 10^{-9}$ & $1.61\times 10^{-9}$ & $7.25 \times 10^{-10}$& $1.67 \times 10^{-10}$ & $5.01 \times 10^{-9}$ & $9.46 \times 10^{-11}$
		\\
			$P_{preempt}^*$ & $2.54 \times 10^{-5}$ & $2.32 \times 10^{-5}$ & $2.20 \times 10^{-5}$ & $2.63 \times 10^{-5}$ & $2.69 \times 10^{-5}$ & $2.05 \times 10^{-5}$ & $2.03 \times 10^{-5}$
		\\
		\hline
		$\lambda_{\mathcal{N}} = 0.5$,    $\mu_{\mathcal{N}}=1$ & & & & & & &\\
		\hline
	Iterations & 1000 & 1003 &1003 & 1003 & 1003 & 1003 & 1003\\
		$\mu_{\mathcal{H}}$ & 0.5 & 0.625 & 0.75 & 0.875 & 1 & 1.125 & 1.25 \\
	$S^*$ & 3 &3 &3 &3& 3& 3& 3\\
	$\lambda_{\mathcal{H}}^*$ & 0.01001 & 0.01001 & 0.01001 & 0.01001 & 0.01001 & 0.01001 & 0.01001\\
		$P_d^*$ & $1.77 \times 10^{-8}$ & $6.12 \times 10^{-9}$ & $2.51 \times 10^{-9}$ & $1.16 \times 10^{-9}$& $5.83 \times 10^{-10}$ & $3.17 \times 10^{-10}$ & $2.70 \times 10^{-10}$
		\\
		$P_{preempt}^*$ & $7.41 \times 10^{-5}$ & $6.94 \times 10^{-5}$ & $9.25 \times 10^{-5}$ & $6.68 \times 10^{-5}$ & $6.41 \times 10^{-5}$ & $6.35 \times 10^{-5}$ & $8.74 \times 10^{-5}$
		\\
		\hline
		$\lambda_{\mathcal{N}} = 0.6$,    $\mu_{\mathcal{N}}=1$ & & & & & & &\\
		\hline
	Iterations & 1000 & 1003 &1003 & 1003 & 1003 & 1003 & 1003\\
		$\mu_{\mathcal{H}}$ & 0.5 & 0.625 & 0.75 & 0.875 & 1 & 1.125 & 1.25 \\
	$S^*$ & 3 &3 &3 &3& 3& 3& 3\\
	$\lambda_{\mathcal{H}}^*$ & 0.01001 & 0.01001 & 0.01001 & 0.01001 & 0.01001 & 0.01001 & 0.01001\\
		$P_d^*$ & $2.36 \times 10^{-8}$ & $8.40 \times 10^{-9}$ & $3.52 \times 10^{-9}$ & $1.66 \times 10^{-9}$& $8.50 \times 10^{-10}$ & $4.68 \times 10^{-10}$ & $2.72 \times 10^{-10}$
		\\
		$P_{preempt}^*$ & $1.76 \times 10^{-4}$ & $1.67 \times 10^{-4}$ & $1.62 \times 10^{-4}$ & $1.59 \times 10^{-4}$ & $1.57 \times 10^{-4}$ & $1.56 \times 10^{-4}$ & $1.55 \times 10^{-4}$\\
		\hline
	\end{tabular}}
	\caption{Optimal values of $S^*$ and $\lambda_{\mathcal{H}}^*$ for different values of $\mu_{\mathcal{H}}$ and $\lambda_{\mathcal{N}}$ by applying SA method}
	\label{tab:my_label2}
\end{table}

\section{Conclusions} \label{section7}
In the leading edge wireless technologies, traffic of different classes, e.g., video, voice, images, data, etc., are assigned different categories of importance,  and consequently their services are effectuated in accordance with an appropriate priority policy.  In cellular networks, the kinds of systems where a higher priority traffic has an advantage in access to service compared to less important ones, are explored through priority queueing models. Therefore, in this study  a $\textrm{\it MMAP[2]/PH[2]/S}$  preemptive repeat priority queueing model with \textrm{$P\!H$} distributed retrial times is investigated. The incoming traffic is classified into handoff calls and new calls where  preemptive repeat priority  is  assigned  to  handoff  calls  over  new  calls. Due to the brief span of inter-retrial times in comparison to service times, a more generalized approach, \textrm{$P\!H$} distributed retrial
times is used so that the performance of the system is not over or under estimated. The analysis of the proposed model is implemented via investigation of steady-state  behaviour of $\textrm{\it L\!D\!Q\!B\!D}$ process.
%The steady-state analysis for the proposed $\textrm{\it L\!D\!Q\!B\!D}$ process has been carried out by  applying matrix analytic algorithm. 
The expressions for the important performance measures have been derived and successfully implemented to demonstrate the influence of various rates  on the performance measures of the system. Further, to prevent frequent termination of services for new calls due to the consideration of the preemptive repeat priority policy, a traffic control  optimization problem has been formulated to estimate the optimal values of parameters $\lambda_{\mathcal{H}}$ and $S$ such that the loss probabilities $P_d$ and $P_{preempt}$ must not exceed some pre-defined threshold values. The  proposed optimization problem has been investigated  by employing \textit{DS}, \textit{PSO} and \textit{SA} methods. The obtained identical numerical results exhibit the reliability of its  optimal solutions. The presented results of the optimization problems can be used in modern wireless cellular networks where the flows of traffic might be essentially heterogeneous with respect to arrival and service processes. In the future, authors propose to extend this model  by using the preemptive resume priority policy for various classes of traffic.

\end{document}